\documentclass[12pt]{article}
\usepackage[english]{babel}
\usepackage[utf8]{inputenc}
\usepackage[T1]{fontenc}
\usepackage{amsmath}
\usepackage{amssymb}
\usepackage{amsfonts}
\usepackage{amsthm}
\usepackage{bm}
\usepackage{bbm}
\usepackage{hyperref}
\usepackage{upgreek}

\newcommand{\eps}{\varepsilon}
\newcommand{\C}{\mathbb{C}}
\newcommand{\N}{\mathbb{N}}

\newcommand{\R}{\mathbb{R}}

\newcommand{\Z}{\mathbb{Z}}
\newcommand{\bB}{\bm{B}}

\newcommand{\boC}{\mathcal{C}}

\newcommand{\boE}{\mathcal{E}}
\newcommand{\boF}{\mathcal{F}}

\newcommand{\boI}{\mathcal{I}}

\newcommand{\boM}{\mathcal{M}}

\newcommand{\boO}{\mathcal{O}}
\newcommand{\boP}{\mathcal{P}}
\newcommand{\boQ}{\mathcal{Q}}
\newcommand{\boR}{\mathcal{R}}

\newcommand{\boU}{\mathcal{U}}
\newcommand{\boV}{\mathcal{V}}

\renewcommand{\Re}{\mathop{{\rm Re}}\nolimits}

\newtheorem{cor}{Corollary}
\newtheorem{lem}{Lemma}
\newtheorem{prop}{Proposition}
\newtheorem{thm}{Theorem}
\newtheorem{extthm}{Theorem}

\theoremstyle{definition}
\newtheorem*{merci}{Acknowledgments}

\newtheorem{step}{Step}

\begin{document}

\title{On the stability of the Ginzburg-Landau vortex}

\author{
\renewcommand{\thefootnote}{\arabic{footnote}}
Philippe Gravejat\footnotemark[1],~ Eliot Pacherie\footnotemark[2]~ and
Didier Smets\footnotemark[3]}
\footnotetext[1]{CY Cergy Paris Universit\'e, Laboratoire de Math\'ematiques AGM, 
F-95302 Cergy-Pontoise Cedex, France. E-mail: {\tt philippe.gravejat@cyu.fr}}
\footnotetext[2]{NYUAD Research Institute, New York University Abu Dhabi, 
PO Box 129188, Abu Dhabi, UAE. E-mail: {\tt ep2699@nyu.edu}}
\footnotetext[3]{Sorbonne Universit\'e, Laboratoire Jacques-Louis Lions, Bo\^ite
Courrier 187, 75252 Paris Cedex 05, France. E-mail: {\tt
didier.smets@sorbonne-universite.fr}}
\maketitle

\abstract{We introduce a functional framework taylored to investigate the 
minimality and stability properties of the Ginzburg-Landau vortex of degree one 
on the whole plane. We prove that a renormalized Ginzburg-Landau energy is 
well-defined in that framework and that the vortex is its unique global 
minimizer up to the invariances by translation and phase shift. Our main result 
is a nonlinear coercivity estimate for the renormalized energy around the 
vortex, from which we can deduce its orbital stability as a solution to the 
Gross-Pitaevskii equation, the natural Hamiltonian evolution equation associated 
to the Ginzburg-Landau energy.}

\section{Introduction}

We are interested in the complex Ginzburg-Landau equation in the plane 
\begin{equation}
\label{eq:GL}
\Delta \Psi + (1 - |\Psi|^2) \Psi = 0. 
\end{equation}
For each $d \in \Z^*$, this equation possesses a well-known 
solution called the vortex of degree $d$ at infinity. It has the equivariant 
form
\begin{equation}
\label{eq:vortex}
V_d(x) = \rho_d(r) e^{i d \theta}, 
\end{equation}
for $x = (r \cos(\theta), r \sin(\theta))$. The profile $\rho_d$ is real-valued, 
increasing, smooth, and satisfies $\rho_d(0) = 0$ and $\rho_d(r) \to 1$ as $r 
\to + \infty$. Equation~\eqref{eq:GL} is invariant by translations and by 
constant phase shifts, so that functions of the form $e^{i \varphi} V_d(\cdot - 
a)$ for arbitrary $a \in \R^2$ and $\varphi \in \R$ are also solutions. It is 
also invariant by complex conjugacy. This is reflected in the equality $V_d = 
\bar{V}_{- d}$. For this reason, we restrict in the sequel our attention to the 
case $d \geq 1$.

Associated to the Ginzburg-Landau equation is the Ginzburg-Landau energy
$$ 
\boE_\text{GL}(\Psi) \equiv \int_{\R^2} e_\text{GL}(\Psi) := \int_{\R^2} \Big( 
\frac{1}{2} |\nabla \Psi|^2 + \frac{1}{4} (1 - |\Psi|^2)^2 \Big).
$$
The vortices $V_d$ do \textit{not} have finite energy. This has long been a source of 
difficulty for their analysis, leading to various strategies based on suitable 
forms of \textit{renormalization}. The overall picture is that $V_1$ possesses some 
minimizing and stability properties, while $V_d$ is unstable for any $d \geq 2$. 
 Over the years, this has been shown in different frameworks, some of which will 
be recalled below, depending on what is precisely meant by minimality and/or 
stability. 

Regarding the stability of $V_1$, a very natural question is that of its 
nonlinear dynamical stability as a stationary solution to the corresponding 
Hamiltonian evolution equation, the Gross-Pitaevskii equation
\begin{equation} 
\label{eq:GP}
i \partial_t \Psi + \Delta \Psi + (1 - |\Psi|^2) \Psi = 0. 
\end{equation}

Our goal in this work is two-fold. First, improving on some of the existing 
variational stability estimates for $V_1$. This will involve nonlinear 
coercivity estimates for (a renormalized version of) $\boE_\text{GL}$ around 
$V_1$. Second, proving the orbital stability of $V_1$ as a solution 
to~\eqref{eq:GP} in a natural energy space. This was an open problem even for 
smooth and compactly supported initial perturbations.

In the remaining part of this introduction, we present the functional framework 
and state our main results. The strategy leading to these results and how they 
relate to earlier works in the literature is the object of the next section.

We introduce the complex Hilbert space
$$
H := \big\{ \Psi \in L^2_\text{loc}(\R^2, \C) \text{ s.t. } \| \Psi \|_{H} < + \infty \big\},
$$
corresponding to the norm
$$
\| \Psi \|_H^2 := \int_{\R^2} \big( |\nabla(\Psi \bar{V}_1)|^2 + (1 - |V_1|^2) 
|\nabla \Psi|^2 \big),
$$
where we recall that $1 - |V_1|^2 = 1 - \rho_1^2 > 0$. We define the energy 
space $E$ as
$$
E := \big\{ \Psi \in H \text{ s.t. } 1 - |\Psi|^2 \in L^2(\R^2) \big\}.
$$
It is a complete metric space for the distance
$$
d_{E}(\Psi_1, \Psi_2) := \| \Psi_1 - \Psi_2 \|_H + \| |\Psi_1|^2 - |\Psi_2|^2 \|_{L^2}.
$$

The vortex $V_1$ belongs to $E$. An important feature of this space is that the 
infinitesimal generators of the invariance groups, i.e. $\partial_{x_1} V_1$ and 
$\partial_{x_2} V_1$ for the translations, and $i V_1$ for the phase shifts, all 
belong to $H$. Actually, although the norm $\| \cdot \|_H$ itself is not 
translation invariant, the spaces $H$ and $E$ are invariant by translations and 
by constant phase shifts. 


Our first result shows that the space $E$ is a natural framework for the 
renormalization of the Ginzburg-Landau energy with respect to $V_1$. 

\begin{prop}
\label{prop:renorm}
The renormalized Ginzburg-Landau energy 
\begin{equation}
\label{def:renorm}
\boE(\Psi) := \lim_{r \to + \infty} \int_{B_r} \big( e_\text{GL}(\Psi) - e_\text{GL}(V_1) \big)
\end{equation}
is well-defined on $E$. Besides, it is invariant by translations and constant phase shifts. 
\end{prop}

It turns out that the condition $\Psi \in E$ actually \textit{encodes} the fact 
that $\Psi$ has a ``degree one'' at infinity, even when the zero set of $\Psi$ 
might be unbounded and the notion of topological degree in the classical meaning 
would make no direct sense. The existence of the limit defining the renormalized 
energy $\boE$ in Proposition~\ref{prop:renorm} is also a consequence of the fact 
that $\Psi \in E$, and it could not be used as a standalone ``definition'' of an 
appropriate functional framework~\footnote{The limit might exist for oscillating 
pure phases, e.g., although such fields would have a well-defined topological 
degree being zero, not one.}. 

Regarding the minimality of $V_1$, we show 

\begin{prop}
\label{prop:minim}
The vortex $V_1$ is the unique global minimizer 
of the energy $\boE$ on $E$, up to translations and constant phase shifts. 
\end{prop}

Our next result, and the core of this work, is a coercivity estimate for the 
renormalized energy in $E$. We denote the orbit under the invariance groups of a 
function $\Psi \in E$ by
$$
\text{Orb}(\Psi) := \big\{ e^{- i \varphi} \Psi(\cdot + a), \varphi \in \R \text{ and } a \in \R^2\}.
$$
Notice that, by construction, $\boE(V_1) = 0$.

\begin{thm}
\label{thm:estimcore}
There exist $\kappa > 0$ and $\rho > 0$ such that
$$
\boE(\Psi) \geq \kappa \, d_E(V_1, \text{Orb}(\Psi))^2,
$$
whenever $\Psi \in E$ satisfies $d_E(V_1, \text{Orb}(\Psi)) < \rho$.
\end{thm}

The distance $d_E(V_1, \text{Orb}(\Psi))$ vanishes if and only if $\Psi \in \text{Orb}(V_1)$. 
Therefore, Theorem~\ref{thm:estimcore} is a (nonlinear) coercivity estimate 
``perpendicular'' to the orbit of $V_1$. 

In this statement, it might have appeared more familiar to use the quantity $d_E(\Psi, 
\text{Orb}(V_1))$ instead of $d_E(V_1, \text{Orb}(\Psi))$, but the former is not 
invariant under translations of $\Psi$, and the norm $\| \cdot \|_H$, which is part of 
the definition of $d_E$, is taylored to measure perturbations around $V_1$. 
Instead, the inequality in Theorem~\ref{thm:estimcore} is invariant by 
translation, since both $\boE(\Psi)$ and $\text{Orb}(\Psi)$ are.

\medskip

Concerning the Gross-Pitaevskii equation~\eqref{eq:GP}, the affine space $\Psi_0 
+ H^1(\R^2, \C)$ is contained in $E$ for arbitrary $\Psi_0 \in E$, and we can 
prove

\begin{prop}
\label{prop:cauchy}
For any $\Psi_0 \in E$, the Gross-Pitaevskii equation~\eqref{eq:GP} possesses a 
unique global solution $t \mapsto \Psi_t \in \boC^0(\R, \Psi_0 + H^1(\R^2, \C))$ 
with initial datum $\Psi_0$. Moreover, the renormalized energy $\boE$ is 
conserved along the flow.
\end{prop}

It is conceivable that the Cauchy problem for~\eqref{eq:GP} is actually globally 
well-posed on $E$, and not just on fibers of the form $\Psi_0 + H^1(\R^2, \C)$. 
We have not settled that question, but it is not necessary to prove our 
stability result. Combining Theorem~\ref{thm:estimcore} and 
Proposition~\ref{prop:cauchy}, we indeed deduce

\begin{thm}
\label{thm:stab}
There exist $\delta > 0$ and $C > 0$ such that, if $\Psi_0 \in E$ and $d_E(V_1, 
\Psi_0) \leq \delta$, then the solution $\Psi_t$ with initial datum $\Psi_0$ of 
Proposition~\ref{prop:cauchy} satisfies
$$
d_E \big( V_1, \text{Orb}(\Psi_t) \big) \leq C \, d_E(V_1, \Psi_0),
$$
for any $t \in \R$. In particular, the vortex $V_1$ is orbitally stable.
\end{thm}

The statement in Theorem~\ref{thm:stab} does not provide any information on the 
location of the solution with respect to the orbit. In the course of the proof, 
we actually construct positions $a(t) \in \R^2$ and phase shifts $\varphi(t) \in 
\R$ such that
$$
d_E \big( V_1, e^{- i \varphi(t)} \Psi_t(\cdot + a(t)) \big) \leq C d_E(V_1, \Psi_0),
$$
for any $t \in \R$. The following proposition provides a first control on the 
evolution of these parameters.

\begin{prop}
\label{prop:evol-param}
There exist $\tau > 0$ and $C > 0$ such 
that, if $\Psi_0 \in E$ and $d_E(V_1, \Psi_0) \leq \tau$, then there exist two 
functions $a \in \boC^1(\R, \R^2)$ and $\varphi \in \boC^1(\R, \R)$ such that
$$ 
d_E \big( V_1, e^{- i \varphi(t)} \Psi_t(\cdot + a(t)) \big) \leq C d_E(V_1, 
\Psi_0),
$$
and
\begin{equation}
\label{eq:est-evol-modul}
\big| a'(t) \big| + \big| \varphi'(t) \big| \leq C d_E(V_1, \Psi_0),
\end{equation}
for any $t \in \R$.
\end{prop}

The question whether $V_1$ is stable as a stationary solution to~\eqref{eq:GP}, 
and not only orbitally stable, is still open. There is no immediate obstruction 
to that stronger form of stability since, although there exist travelling waves 
of~\eqref{eq:GP} with arbitrarily small speed (see~\cite{BethSau1, ChirPac1}), 
they are not small perturbations of the vortex but instead perturbations of a 
vortex-antivortex pair, and have finite Ginzburg-Landau energy.

Asymptotic stability of $V_1$ (or maybe only of its orbit) is also an open 
question. Although~\eqref{eq:GP} is an Hamiltonian equation, asymptotic 
stability could hold for a topology in which the renormalized energy is not 
continuous (e.g. through dispersion at infinity). In one space dimension, it was 
proved in~\cite{GravSme1} that the black soliton, the 1d equivalent of $V_1$, is 
asymptotically stable in an orbital sense.

Finally, we mention that the dynamics of~\eqref{eq:GP} for sequences of initial 
data that converge to $V_1$ or more generally to suitable combinations of well 
separated vortices of degrees $\pm 1$ has already been studied on finite time 
intervals (see e.g.~\cite{CollJer1, LinXin1, BetJeSm1, JerrSpi2} and the references 
therein). In particular, it is known that a related notion of modulation 
parameters asymptotically obey a limit point-vortex ODE, but only on finite time 
intervals. These results are based on different kind of rigidity estimates for $V_1$
, where the closeness is measured through the Jacobian (see \cite{JerrSpi1} for
the sharpest known statement in that direction). 
It is tempting to investigate whether the results of the present 
paper could be used to extend the description of this dynamics to longer time 
scales.

\medskip

In the sequel, we use all along the following notation. We set $x^\perp = (- 
x_2, x_1)$ for any $x \in \R^2$. We also use the notation $B_r$ for the open 
ball of $\R^2$ with center $(0, 0)$ and radius $r > 0$. Finally, we define the 
scalar product of two complex numbers $z_1 = a_1 + i b_1$ and $z_2 = a_2 + i 
b_2$ as $\langle z_1, z_2 \rangle_\C = \Re (z_1 \overline{z_2}) = a_1 a_2 + b_1 
b_2$.

\section{Strategy for the proofs}

In this section, we present the key ingredients that are needed in the proofs
of the results stated in the introduction. This is also the occasion to
discuss how our arguments relate to earlier works in the literature. The
detailed proofs or their completions is postponed to the subsequent sections.

\subsection{Concerning renormalization and Proposition~\ref{prop:renorm}}

We start by motivating the introduction of the Hilbert space $H$. At the same 
time, we show how its definition implies the existence of the limit 
in~\eqref{def:renorm} providing the renormalized Ginzburg-Landau energy $\boE$.

To do so, we first recall the well-known fact that the divergence of the 
Ginzburg-Landau energy of $V_1$ is \textit{only} due to a too slow decay of the 
gradient of its phase at infinity, because of its non-trivial topological 
degree. Taking the gradient of~\eqref{eq:vortex}, we obtain 
\begin{equation}
\label{eq:nabla-V1}
\nabla V_1(x) = e^{i\theta} \Big( \rho_1'(|x|) \frac{x}{|x|} + i \rho_1(|x|) 
\frac{x^\perp}{|x|^2} \Big), 
\end{equation}
and therefore
$$
\frac{1}{2} \int_{B_r} |\nabla V_1|^2 = \pi \int_0^r \Big( \rho_1'(r)^2 +
\frac{\rho_1(r)^2}{r^2} \Big) \, r \, dr.
$$
From the asymptotic properties of $\rho_1$, which we have recalled in 
Lemma~\ref{lem:propprofil} in Appendix~\ref{sec:V-1}, it follows that 
$$
\pi \int_0^r \frac{\rho_1(r)^2}{r^2} \, r \, dr = \pi \log(r) + \boO(1),
$$
as $r \to + \infty$, while
\begin{equation}
\label{eq:rho1good}
\int_0^{+ \infty} \Big( \frac{\rho_1'(r)^2}{2} + \frac{(1 - \rho_1(r)^2)^2}{4} 
\Big) \, r \, dr < + \infty.
\end{equation}
The definition of the space $H$ through the norm
$$
\| \Psi \|_H^2 := \int_{\R^2} \big( |\nabla(\Psi \bar{V}_1)|^2 + (1 - |V_1|^2) 
|\nabla \Psi|^2 \big),
$$
is then related to our requirement that the infinitesimal generators 
$\partial_{x_1} V_1,$ $\partial_{x_2} V_1$ and $iV_1$ of the invariance groups 
should all belong to $H$ in order to eventually derive \textit{optimal} 
coercivity estimates. It follows from the properties listed in 
Lemma~\ref{lem:propprofil} that the derivatives $\partial_{x_1} V_1$ and 
$\partial_{x_2} V_1$ have gradients in $L^2(\R^2)$. As we have just recalled, 
the same does \textit{not} hold for the function $i V_1$. On the other hand, if 
$\Psi = i V_1$ then $\Psi \bar{V}_1 = i |V_1|^2$, and the latter has a gradient 
in $L^2(\R^2)$, in view of~\eqref{eq:rho1good} and the boundedness of $|V_1|$.

From the Leibniz rule, we compute
\begin{equation}
\label{eq:exp-nabla-V1-Psi}
|\nabla(\Psi \bar{V}_1)|^2 = |\nabla \Psi|^2 |V_1|^2 + |\Psi|^2 |\nabla V_1|^2 + 
2 \langle \nabla \Psi \bar{V}_1, \Psi \nabla \bar{V}_1\rangle_\C, 
\end{equation}
where, here as in the sequel, we denote by $\langle \cdot, \cdot \rangle_\C$ the 
scalar product of $\C \simeq \R^2$. With the definition~\eqref{def:renorm} of 
the renormalized Ginzburg-Landau energy in mind, we write
$$
|\nabla \Psi|^2 - |\nabla V_1|^2 = |\nabla(\Psi \bar{V}_1)|^2 + (1 - |V_1|^2)| 
\nabla \Psi|^2 - (1 - |\Psi|^2) |\nabla V_1|^2 - 2 \langle \nabla(\Psi 
\bar{V}_1), \Psi \nabla \bar{V}_1\rangle_\C.
$$
Note that the first two terms in the right-hand side of the previous pointwise 
identity are precisely those entering in the definition of the norm in $H$, and 
hence by construction they are integrable over $\R^2$ when $\Psi \in H$. The 
third term is also integrable when in addition $\Psi \in E$. This follows from 
the Cauchy-Schwarz inequality since $1 - |\Psi|^2 \in L^2(\R^2)$ for $\Psi \in 
E$, and on the other hand $\nabla V_1 \in L^4(\R^2)$. In order to prove the 
existence of the limit in~\eqref{def:renorm}, it therefore only remains to prove 
the existence of the limit
\begin{equation}
\label{eq:mom}
\lim_{r \to +\infty} \int_{B_r} \langle \nabla(\Psi \bar{V}_1),
\Psi \nabla \bar{V}_1 \rangle_\C,
\end{equation}
whenever $\Psi \in E$. This requires the identification of some cancellation 
phenomenon. For this and also later purposes, we split the integral 
in~\eqref{eq:mom} into a local part and a part at infinity. More precisely, here 
and throughout the paper, we fix a smooth, radial, radially non-increasing 
cut-off function $0 \leq \chi \leq 1$ such that $\chi \equiv 1$ in $B_1$ and 
$\chi$ is supported in $B_2$. For arbitrary $R > 0$, we set $\chi_R(x) := 
\chi(x/R)$. Having in mind~\eqref{eq:nabla-V1}, we decompose $\langle 
\nabla(\Psi \bar{V}_1), \Psi \nabla \bar{V}_1 \rangle_\C$ as
$$
\Big\langle \nabla(\Psi \bar{V}_1), \Big( \nabla \bar{V}_1 + i (1 - \chi_R)^2 
\frac{x^\perp}{|x|^2} \bar{V}_1 \Big) \Psi \Big\rangle_\C - 
\frac{x^\perp}{|x|^2} (1 - \chi_R)^2 \langle i \Psi \bar{V}_1, \nabla(\Psi 
\bar{V}_1) \rangle_\C.
$$

The first term in the previous decomposition is integrable on $\R^2$. As a 
matter of fact, it follows from~\eqref{eq:nabla-V1} that
$$
\Big( \nabla \bar{V}_1 + i (1 - \chi_R)^2 \frac{x^\perp}{|x|^2} \bar{V}_1 \Big) 
\Psi = e^{i \theta} \Big( \rho_1' \frac{x}{|x|} + i \chi_R (2 - \chi_R) \rho_1 
\frac{x^\perp}{|x|^2} \Big) \Psi,
$$
and the latter quantity belongs to $L^2(\R^2)$. Indeed, we first infer from 
Lemma~~\ref{lem:propprofil} again that $\rho_1'(r) \sim 1/r^3$, as $r \to + 
\infty$. Second, the function $\chi_R$ has compact support, while the function 
$\rho_1 x^\perp/|x|^2$ is bounded. The previous claim finally results from the 
fact that $\Psi$ is in the space $H$, which continuously embeds into the 
weighted space $L_{- s}^2(\R^2)$ for arbitrary $s > 1$ (as shown in 
Lemma~\ref{lem:comp-emb}). Here, we only need the case $s = 3$. 

At this stage, we have reduced the existence of the renormalized energy to the 
next key lemma, where we additionally gain some smallness estimate in the 
vicinity of $V_1$. Note that this gain will be important for our later 
perturbation analysis.

\begin{lem}
\label{lem:errR}
Let $\Psi \in E$. The quantity
\begin{equation}
\label{def:errR}
P_R(\Psi) := \lim_{r \to + \infty} 2 \int_{B_r} (1 - \chi_R)^2 
\frac{x^\perp}{|x|^2} \cdot \langle i \Psi \bar{V}_1, \nabla(\Psi \bar{V}_1) 
\rangle_\C
\end{equation}
is well-defined for any $R > 0$. Moreover, there exist universal constants 
$\delta > 0$, $\Lambda > 0$ and $K \geq 1$ such that
\begin{equation}
\label{eq:smallerrR}
|P_R(\Psi)| \leq \frac{K}{R} \, d_E(\Psi, V_1)^2,
\end{equation}
provided that $d_E(\Psi, V_1) \leq \delta$ and $R \geq \Lambda$.
\end{lem}

The proof of the first statement in Lemma~\ref{lem:errR} is concise, and we 
present it here next. For $\Psi \in E$, it is shown in Lemma~\ref{lem:finiGL} 
that $\Psi \bar{V}_1$ has finite Ginzburg-Landau energy. Therefore, it follows 
from a result of P.~G\'erard~\cite{Gerard2} that $\Psi \bar{V}_1$ may be 
decomposed as 
\begin{equation}
\label{eq:PG}
\Psi\bar{V}_1 = e^{i\varphi} + w,
\end{equation}
for some real-valued function $\varphi$ such that $\nabla \varphi \in L^2(\R^2)$ 
and some complex-valued function $w \in H^1(\R^2)$. We then split the scalar 
product in~\eqref{def:errR} as
$$
\langle i \Psi \bar{V}_1, \nabla(\Psi \bar{V}_1) \rangle_\C = \langle 1, \nabla 
\varphi \rangle_\C + \langle i w, \nabla w \rangle_\C + \langle w, \nabla 
\varphi e^{i \varphi} \rangle_\C + \langle i e^{i \varphi}, \nabla w \rangle_\C,
$$
and we treat separately each of the four corresponding terms. Given any $r > 0$, 
we first write
$$
\int_{B_r} (1 - \chi_R)^2 \frac{x^\perp}{|x|^2} \cdot \langle 1, \nabla \varphi 
\rangle_\C = \int_{B_r} \text{div}\big( (1 - \chi_R)^2 \frac{x^\perp}{|x|^2} 
\langle 1, \varphi \rangle_\C \big) = 0.
$$
Here, we have used the fact that $\nabla \chi_R(x) \cdot x^\perp = 0$ pointwise 
since $\chi_R$ is radial, and also that the flux of $x^\perp$ through $\partial 
B_r$ vanishes pointwise. Next, the second and third terms are vector fields, 
which are integrable over $\R^2$, since $w \in H^1(\R^2)$ and $\nabla \varphi 
\in L^2(\R^2)$. Therefore, their integral against the bounded vector field $(1 - 
\chi_R)^2 x^\perp/|x|^2$ is well-defined on $\R^2$. Finally, we use as above 
that
$$
\int_{B_r} \text{div} \big( (1 - \chi_R)^2\frac{x^\perp}{|x|^2} \cdot \langle i 
e^{i \varphi}, w \rangle_\C \big) = 0,
$$
in order to obtain the identity for the fourth term
$$
\int_{B_r} (1 - \chi_R)^2 \frac{x^\perp}{|x|^2} \cdot \langle i e^{i \varphi}, 
\nabla w \rangle_\C = \int_{B_r} (1 - \chi_R)^2 \frac{x^\perp}{|x|^2} \cdot 
\langle \nabla \varphi e^{i \varphi}, w \rangle_\C.
$$
The last integral now has a well-defined limit when $r \to + \infty$, since both 
$w$ and $\nabla \varphi$ belong to $L^2(\R^2)$. This completes the proof of the 
existence of $P_R(\Psi)$ in Lemma~\ref{lem:errR}, and therefore also of 
$\boE(\Psi)$ in Proposition~\ref{prop:renorm}. We refer to 
Section~\ref{sect:P1remaining} below for the proof of the invariance with 
respect to translations and phase shifts of this latter quantity. Note also 
that, with the help of decomposition~\eqref{eq:PG}, we have obtained the formula
\begin{equation}
\label{eq:equiverrR}
P_R(\Psi) = 2 \int_{\R^2} (1 - \chi_R)^2 \frac{x^\perp}{|x|^2} \cdot \langle i 
w, \nabla w+ 2 i \nabla \varphi e^{i \varphi} \rangle_\C.
\end{equation}

The proof of~\eqref{eq:smallerrR} follows similar lines, the main difference 
being that decomposition~\eqref{eq:PG} needs to be adapted to a perturbative 
setting. This is done in Lemma~\ref{lem:PGalt} of Section 
\ref{sect:P1remaining}. Note that the integrand in~\eqref{def:errR} identically 
vanishes for $\Psi = V_1$, since $\langle i |V_1|^2, \nabla(|V_1|^2) \rangle_\C 
= 0_{\R^2}$ due to the real-valued nature of the function $|V_1|$. This may 
serve as an intuition to why~\eqref{eq:smallerrR} actually holds.

For later reference, we also make here explicit the decomposition of 
$\boE(\Psi)$ which we have obtained so far, namely
\begin{equation}
\label{eq:decompboE}
\begin{split}
\boE(\Psi) = & \frac{1}{2} \| \Psi \|_H^2 - \int_{\R^2} |\nabla V_1|^2 (1 - 
|\Psi|^2)\\
& - \int_{\R^2} \Big\langle \nabla(\Psi \bar{V}_1), \Big( \nabla 
\bar{V}_1 + i (1 - \chi_R)^2 \frac{x^\perp}{|x|^2} \bar{V}_1 \Big) \Psi 
\Big\rangle_\C + \frac{1}{2} P_R(\Psi)\\
& + \int_{\R^2} \frac{1}{4} \big( (1 - |\Psi|^2)^2 - (1 - |V_1|^2)^2 \big).
\end{split}
\end{equation}

\subsection{Concerning minimality and Proposition~\ref{prop:minim}}

A solution $\Psi$ to~\eqref{eq:GL} is called a locally minimizing solution if
$$
\boE_\text{GL}(\Psi + \varepsilon, B_R) \geq \boE_\text{GL}(\Psi, B_R),
$$
for any $R >0$ and any $\varepsilon \in H_0^1(B_R, \C)$. Here, we have set
$$
\boE_\text{GL}(\Psi, B_R) := \int_{B_R} e_\text{GL}(\Psi).
$$
The next characterization was obtained by P. 
Mironescu~\cite{Mirones1}~\footnote{The fact that the vortex solution $V_1$ is 
locally minimizing is not explicitly stated in~\cite{Mirones1}, but it follows 
from properties listed in there, in particular Corollaire 2 and the remark 
following it.}.

\begin{extthm}[\cite{Mirones1}] 
\label{thm:MS}
The vortex solution $V_1$ is a locally minimizing solution to~\eqref{eq:GL}. 
Moreover, it is the only non-constant locally minimizing solution, up to 
translations and constant phase shifts.
\end{extthm}

Since any minimizer of $\boE$ in $E$ is necessarily also a locally minimizing 
solution to~\eqref{eq:GL}, and since $\boE(V_1) = 0$ by construction, the proof 
of Proposition~\ref{prop:minim} reduces to show that $\boE$ is non-negative. For 
that purpose, we shall appeal, after suitable rescalings, to results (for 
example~\cite[Corollaire 2]{Mirones1}) regarding the asymptotics of 
Ginzburg-Landau minimizers on a fixed bounded domain with fixed boundary data. 
The reduction to the latter case from our framework requires some elementary 
surgery on the boundary of large balls, the necessary details of which are 
presented in Section~\ref{sect:P1remaining}.

\subsection{Concerning coercivity and Theorem~\ref{thm:estimcore}}

A quantitative stability estimate with respect to compactly supported 
perturbations of $V_1$ was also obtained by P. Mironescu in~\cite{Mirones2}. For 
that purpose, he decomposed
\begin{equation}
\label{eq:decompmiro}
\boE_\text{GL}(V_1 + \varepsilon, B_R) = \boE_\text{GL}(V_1, B_R) + \frac{1}{2} 
B(\varepsilon) + \boO \big( \| \varepsilon \|_{H_0^1(B_R)}^3 \big),
\end{equation}
for any function $\varepsilon \in H_0^1(B_R)$. Here, $B$ is the real quadratic 
form given by
$$
B(\varepsilon) := \int_{\R^2} \big( |\nabla \varepsilon|^2 - (1 - |V_1|^2) 
|\varepsilon|^2 + 2 \langle V_1, \varepsilon \rangle_\C^2 \big).
$$

\begin{extthm}[\cite{Mirones2}]
For any $R > 0$, there exists $\kappa_R > 0$ such that 
$$
B(\varepsilon) \geq \kappa_R \, \| \varepsilon \|_{H_0^1(B_R)}^2,
$$
for any $\varepsilon \in H_0^1(B_R)$.
\end{extthm}

The fact that the invariance by translation and by phase shift is not reflected 
in the previous coercivity estimate is due to the restriction $\varepsilon \in 
H_0^1(B_R)$, which prevents those groups to act. In turn, this can be used to 
show that necessarily $\kappa_R \to 0$ as $R \to + \infty$.

In order to derive a stability estimate without restricting to compactly 
supported perturbations, M. del Pino, P. Felmer and M. Kowalczyk~\cite{dePiFeK1} 
considered an Hilbert space $H_B$ naturally associated to the decomposition 
\eqref{eq:decompmiro} (see also~\cite{OvchSig1} for previous approach in 
the space $L^2(\R^2)$). 
This space was defined from the norm~\footnote{Note the 
sign change in the middle term with respect to the quadratic form $B$.}
$$
\| \varepsilon \|_{H_B}^2 := \int_{\R^2} \big( |\nabla \varepsilon|^2 + (1 - 
|V_1|^2) |\varepsilon|^2 + 2 \langle V_1, \varepsilon \rangle_\C^2 \big).
$$
One can note that $H_B$ is a strict subspace of $H$. With the definition of this 
space at hand, they obtained

\begin{extthm}[\cite{dePiFeK1}]
\label{thm:dfk}
The real quadratic form $B$ is positive semi-definite on $H_B$, and its kernel 
coincides with the real vector space spanned by $\partial_{x_1}V_1$ and 
$\partial_{x_2}V_1$.
\end{extthm}

The invariance by translation of the Ginzburg-Landau energy is reflected in 
Theorem~\ref{thm:dfk} in the fact that $\partial_{x_1} V_1$ and $\partial_{x_2} 
V_1$ belong to the kernel of $B$. The invariance by phase shifts cannot be 
accounted for by working in the space $H_B$ though, since $iV_1 \notin H_B$. 
This is the reason why the latter is not present in the kernel of $B$ in $H_B$. 

It is possible to extend Theorem~\ref{thm:dfk} with some quantitative coercivity 
estimates for $B$ under suitable orthogonality conditions of $\varepsilon$ with 
respect to $\partial_{x_1} V_1$, $\partial_{x_2} V_1$ and $iV_1$ (see 
e.g.~\cite[Proposition 1.3]{ChirPac2} and also the Fredholm alternative 
in~\cite[Theorem 2]{dePiFeK1}). These estimates however do not allow to control 
the nonlinear terms arising in the expansion of the renormalized Ginzburg-Landau 
energy, and it does not seem possible to derive nonlinear stability of $V_1$ 
based (exclusively) on the linear analysis of $B$.

At this stage, it is worth comparing the quadratic form $B$ with our previous 
decomposition~\eqref{eq:decompboE} of $\boE$. For that purpose, we first write 
$\Psi = V_1 + \varepsilon$ in~\eqref{eq:decompboE}. Using the fact that $V_1$ is 
a solution to~\eqref{eq:GL}, we obtain

\begin{lem}
\label{lem:decomp}
For $\Psi = V_1 + \varepsilon \in E$, we have
\begin{equation}
\label{eq:decompboE2}
\boE(V_1 + \varepsilon) = \frac{1}{2} \boQ_R(\varepsilon) + \frac{1}{2} 
\boP_R(\varepsilon) + \frac{1}{4}\int_{\R^2} \eta_\varepsilon^2.
\end{equation}
In this identity, $\boQ_R(\varepsilon)$ is the quadratic form on $H$ given by
\begin{equation}
\label{eq:defboQ}
\begin{split}
\boQ_R(\varepsilon) := & \big\| \varepsilon \big\|_H^2 - \int_{\R^2} \big( 1 
-|V_1|^2 - |\nabla V_1|^2 \big) |\varepsilon|^2 \\
& - 2 \int_{\R^2} \Big\langle \nabla (\varepsilon \bar{V}_1), \Big( \nabla 
\bar{V}_1 + i \frac{x^\perp}{|x|^2} (1 - \chi_R)^2 \bar{V}_1 \Big) \varepsilon 
\Big\rangle_\C,
\end{split}
\end{equation}
$\boP_R(\varepsilon) := P_R(V_1 + \varepsilon)$, where $P_R$ is defined in 
Lemma~\ref{lem:errR}, and
\begin{equation}
\label{eq:defetaeps}
\eta_\varepsilon := \big( 1 - |V_1 + \varepsilon|^2 \big) - \big( 1 - |V_1|^2 
\big) = - 2\langle \varepsilon, V_1 \rangle_\C - |\varepsilon|^2.
\end{equation}
\end{lem}

An important feature concerning the decomposition in~\eqref{eq:decompboE2} is that 
$\eta_\varepsilon^2$ being a square, it is (pointwise) non-negative. If we 
develop $\eta_\varepsilon^2$ according to definition~\eqref{eq:defetaeps}, we 
obtain the identity
$$
B(\varepsilon) = \boQ_R(\varepsilon) + \boP_R(\varepsilon) + 2 \int_{\R^2} 
\langle \varepsilon, V_1\rangle_\C^2,
$$
but the latter only makes sense provided that $\varepsilon \in H_B$, because of 
the third term. Besides, as we have already mentioned, the coercivity properties 
of $B$ are insufficient to derive the nonlinear stability of $V_1$.

We modify the previous strategy in two ways. First, we only develop the square 
of $\eta_\varepsilon$ according to~\eqref{eq:defetaeps} locally in space. More 
precisely, we write
$$
\frac{1}{4} \int_{\R^2} \eta_\varepsilon^2 = \frac{1}{4} \int_{\R^2} (1 - 
\chi_R^2) \eta_\varepsilon^2 + \frac{1}{4} \int_{\R^2} \chi_R^2 
\eta_\varepsilon^2 = N_R(\varepsilon) + \boI_R(\varepsilon).
$$
Here, the term 
$$
\boI_R(\varepsilon) := \int_{\R^2} \chi_R^2 \langle \varepsilon, V_1 
\rangle_\C^2
$$
is now a real quadratic form well-defined on $H$, and the nonlinear quantity 
$$
N_R(\varepsilon) := \frac{1}{4} \int_{\R^2} (1 - \chi_R^2) \eta_\varepsilon^2 + 
\frac{1}{4} \int_{\R^2} \chi_R^2 \big( |\varepsilon|^4 + 4 \langle \varepsilon, 
V_1 \rangle_\C |\varepsilon|^2 \big)
$$
contains only non-negative terms except possibly the cubic one in $\varepsilon$, 
which however is localized in $B_{2 R}$.

Second, although the quantity $\boP_R(\varepsilon)$ is quadratic in 
$\varepsilon$, we keep it out of our linear analysis. More precisely, we finally 
write
\begin{equation}
\label{eq:decompgood}
\boE(V_1 + \varepsilon) = \frac{1}{2} Q_R(\varepsilon) + N_R(\varepsilon) + 
\frac{1}{2} \boP_R(\varepsilon),
\end{equation}
where $Q_R := \boQ_R + 2 \boI_R$ is the real quadratic form on $H$ for which we 
shall prove a coercivity estimate. The quadratic form $\boP_R$ will eventually 
be controlled using both $Q_R$ and the nonlinear term $N_R$. A careful analysis 
shows that it cannot be included in $Q_R$, since it would otherwise induce an 
infinite number of negative directions.

The coercivity of $Q_R$ is given by

\begin{prop}
\label{prop:quad}
There exist universal constants $\kappa_0 > 0$ and $N_0 > 0$ such that, given 
any $\varepsilon \in H$ verifying the orthogonality conditions
\begin{equation}
\label{eq:orthocond}
\int_{\R^2} \chi \, \langle \varepsilon, \partial_x V_1 \rangle_\C = \int_{\R^2} 
\chi \, \langle \varepsilon, \partial_y V_1 \rangle_\C = \int_{\R^2} \chi \, 
\langle \varepsilon, i V_1\rangle_\C = 0,
\end{equation}
and any $R_0 \geq 1,$ there exists $R_0 \leq R \leq 2^{N_0}R_0$ such that 
$$
Q_R(\varepsilon) \geq \kappa_0 \bigg( \| \varepsilon \|_H^2 + \int_{\R^2} 
\chi_R^2 \langle \varepsilon, V_1 \rangle_\C^2 \bigg).
$$
\end{prop}

The proof of Proposition~\ref{prop:quad} is presented in Section 
\ref{sect:quad}. In a few words, the idea is to decompose the perturbation 
$\varepsilon$ into a local (compactly supported) part, for which the analysis is 
very similar to the one in~\cite{dePiFeK1}, and a second part at infinity, for 
which the exact form of our decomposition, and in particular the fact that 
$\boP_R$ has been left aside, plays a crucial role. The necessary glue between 
the two analysis explains the fact that the cut-off location $R$ is perturbation 
dependent (in a range with universal extent, though).

Observe that the local orthogonality conditions in~\eqref{eq:orthocond} reflect 
the geometric invariances by translation and phase shift of the Ginzburg-Landau 
energy. In the functional framework of the Hilbert space $H$, the $L^2$-scalar 
products between an arbitrary perturbation $\varepsilon$ and the functions 
$\partial_x V_1$, $\partial_y V_1$, respectively $i V1$, do not necessarily make 
sense. This explains the introduction of the smooth radial cut-off function 
$\chi$ in the three integrals of~\eqref{eq:orthocond}. Note that the choice of 
this special function is somewhat arbitrary. One can check that it can be 
replaced by any smooth, non-negative and non-increasing function, identically 
equal to $1$ on $[0, 1]$, and compactly supported.

Concerning the nonlinear term $N_R$, it is straightforward to derive 

\begin{lem}
\label{LN}
For any $R > 0$, $0 < \kappa < 1$, and $\varepsilon \in H$, we have
$$
N_R (\varepsilon) + \kappa \int_{\R^2} \chi_R^2 \langle \varepsilon , V_1 
\rangle_\C^2 \geq \frac{\kappa}{4} \int_{\R^2} \eta_\varepsilon^2 - \| 
\varepsilon \|_{L^3(B_{2 R})}^3.
$$
\end{lem}

Recalling that $\eta_\varepsilon = - 2 \langle \varepsilon, V_1 \rangle_\C - 
|\varepsilon |^2$, we indeed compute that
\begin{align*}
\frac{1}{4} \eta_\varepsilon^2 - (1 - \kappa) \chi_R^2 \langle \varepsilon, V_1 
\rangle_\C^2 \geq & \frac{\kappa}{4} \eta_\varepsilon^2 + (1 - \kappa) \chi_R^2 
\Big( \frac{1}{4} \eta_\varepsilon^2 - \langle \varepsilon, V_1 \rangle_\C^2 
\Big)\\
\geq & \frac{\kappa}{4} \eta_\varepsilon^2 + (1 - \kappa) \chi_R^2 \langle 
\varepsilon, V_1 \rangle_\C |\varepsilon|^2,
\end{align*}
and the conclusion follows after integration on $\R^2$, using that $|V_1| \leq 
1$ and $|1 - \kappa| \leq 1$.

We are now in position to present the detail leading to a nonlinear coercivity 
estimate for $\boE$ around $V_1$ by combining the results in 
Proposition~\ref{prop:quad} for $Q_R$, Lemma~\ref{LN} for $N_R$, and the second 
statement in Lemma~\ref{lem:errR} for $\boP_R$.

Recall the decomposition in~\eqref{eq:decompgood}. In 
Proposition~\ref{prop:quad}, we fix the value of $R_0$ as
$$
R_0 = \max \Big\{ \Lambda, \frac{16 K}{\kappa_0} \Big\},
$$
where $\Lambda \geq 1$ and $K > 0$ are the universal constants provided by 
Lemma~\ref{lem:errR}, and $\kappa_0 > 0$ is the universal constant provided by 
Proposition~\ref{prop:quad}. Under the orthogonality conditions 
in~\eqref{eq:orthocond}, we can find $R \in [R_0, 2^{N_0} R_0] $ such that
\begin{equation}
\label{conclcoerc}
\frac{1}{2} Q_R(\varepsilon) \geq \frac{\kappa_0}{2} \bigg( \|\varepsilon\|_H^2 
+ \int_{\R^2} \chi_R^2 \langle \varepsilon, V_1 \rangle_\C^2 \bigg).
\end{equation}
By Lemma~\ref{lem:errR} and our choice of $R_0$, we also obtain
$$
\Big| \frac{1}{2} \boP_R(\varepsilon) \Big| \leq \frac{K}{2 R} d_E(V_1 + 
\varepsilon, V_1)^2 \leq \frac{\kappa_0}{16}\Big( \| \varepsilon \|_H^2 + \| 
\eta_\varepsilon \|_{L^2}^2 \Big).
$$
Finally, it follows from Lemma~\ref{LN} with $\kappa = \frac{\kappa_0}{2}$ that
$$
N_R (\varepsilon) + \frac{\kappa_0}{2} \int_{\R^2} \chi_R^2 \langle \varepsilon, 
V_1 \rangle_\C^2 \geq \frac{\kappa_0}{8} \| \eta_\varepsilon \|_{L^2}^2 - 
\| \varepsilon \|_{L^3(B_{2 R})}^3.
$$
After summation, the previous three inequalities and~\eqref{eq:decompgood} yield
$$
\boE(V_1 + \varepsilon) \geq \frac{\kappa_0}{16} \Big( \| \varepsilon \|_H^2 + 
\| \eta_\varepsilon \|_{L^2}^2) - \| \varepsilon \|_{L^3(B_{2 R})}^3.
$$
Now, remark that $2 R \leq R_1 := 2^{N_0 +1} R_0$, which is a universal 
constant. Therefore, we derive, in particular from Lemma~\ref{lem:comp-emb}, 
that
$$
\| \varepsilon \|_{L^3(B_{2 R})}^3 \leq \| \varepsilon \|_{L^3(B_{R_1})}^3 \leq 
K_1 \| \varepsilon \|_H^3,
$$
for some further universal constant $K_1 > 0$. If $\| \varepsilon \|_H$ 
satisfies the smallness condition $K_1 \| \varepsilon \|_H \leq \kappa_0/32$, we 
may absorb the remaining cubic term in the quadratic one. Therefore, we have 
proved

\begin{prop}
\label{prop:Ecoerc}
There exist $\kappa > 0$ and $\delta > 0$ such that, given any function $\Psi = 
V_1 + \varepsilon \in E$ such that
$$
\| \varepsilon \|_H + \| \eta_\varepsilon \|_{L^2 (\R^2)} \leq \delta,
$$
and $\varepsilon$ satisfies the three orthogonality conditions 
in~\eqref{eq:orthocond}, we have
$$
\boE (V_1 + \varepsilon) \geq \kappa \big( \| \varepsilon \|_H^2 + \| 
\eta_\varepsilon \|_{L^2}^2 \big).
$$
\end{prop}

Deriving Theorem~\ref{thm:estimcore} from Proposition~\ref{prop:Ecoerc} is then 
mostly a matter of fixing appropriately the orthogonality conditions 
in~\eqref{eq:orthocond}. The strategy to achieve this is classical (see 
e.g.~\cite{Weinste1, Weinste2}). It relies on the introduction of modulation 
parameters corresponding to the geometric invariances. Consider the 
neighbourhoods in $H$ of the orbit of $V_1$ defined as
$$
\boV(\alpha) := \Big\{ \Psi \in H\text{ s.t. }\inf_{(a, \varphi) \in \R^2 \times 
\R} \big\| e^{- i \varphi} \Psi(\cdot + a) - V_1 \big\|_H < \alpha \Big\},
$$
for any $\alpha > 0$. Given a function $\Psi \in \boV(\alpha)$, we decompose it 
as $e^{- i \varphi} \Psi(\cdot + a) = V_1 + \varepsilon$ for $a \in \R^2$ and 
$\varphi \in \R$. When $\alpha$ is small enough, we can choose the modulation 
parameters $a$ and $\varphi$ so as to satisfy the orthogonality conditions 
in~\eqref{eq:orthocond}. More precisely, we shall prove in 
Section~\ref{sec:proof-modul} below

\begin{prop}
\label{prop:modul-param}
There exist $\alpha > 0$ and $A > 0$, and two functions $a \in 
\boC^1(\boV(\alpha), \R^2)$ and $\varphi \in \boC^1(\boV(\alpha), \R / 2 \pi 
\Z)$, such that any function $\Psi \in \boV(\alpha)$ can be written as
\begin{equation}
\label{def:eps}
\Psi = e^{i \varphi(\Psi)} \Big( V_1 \big( \cdot - a(\Psi) \big) + \varepsilon 
\big( \cdot - a(\Psi) \big) \Big),
\end{equation}
where $\varepsilon$ satisfies the orthogonality conditions 
in~\eqref{eq:orthocond}. Moreover, given parameters $(b, \vartheta) \in \R^2 
\times \R$ such that
$$
\big\| e^{- i \vartheta} \Psi(\cdot + b) - V_1 \big\|_H < \alpha,
$$
we have the estimate
\begin{equation}
\label{eq:est-modul}
\big\| \varepsilon \big\|_H + \big| a(\Psi) - b \big| + \big| e^{i 
\varphi(\Psi)} - e^{i \vartheta} \big| \leq A \, \big\| e^{- i \vartheta} 
\Psi(\cdot + b) - V_1 \big\|_H.
\end{equation}
\end{prop}

We may now complete the 

\begin{proof}[Proof of Theorem~\ref{thm:estimcore}]
Let $\Psi \in E$ be such that $d := d_E(V_1, \text{Orb}(\Psi)) < \rho$, where $0 < \rho 
\leq 1$ is a universal constant, which will be fixed in the course of the proof. 
By definition of the distance to the orbit, there exist $\vartheta \in \R$ and 
$b \in \R^2$ such that ~\footnote{We could show that the distance to the orbit 
is actually achieved but we will not need and spare that argument.}
\begin{equation}
\label{eq:almostdist}
\frac{d}{2} \leq d_E \big( V_1, e^{- i \vartheta} \Psi(\cdot + b) \big) \leq d. 
\end{equation}
In particular, we have $\| V_1 - e^{- i \vartheta} \Psi(\cdot + b) \|_H \leq d < 
\rho$. We require that $\rho \leq \alpha$, where $\alpha > 0$ is the constant 
provided by Proposition~\ref{prop:modul-param}. Applying this proposition, we 
derive the existence of $a \in \R^2$ and $\varphi \in \R$ such that $\varepsilon 
:= e^{- i \varphi} \Psi(\cdot +a) - V_1$ satisfies the orthogonality conditions 
in~\eqref{eq:orthocond}, and besides,
$$
\| \varepsilon \|_H + | a - b | \leq A d. 
$$
We next estimate
\begin{align*}
\big\| \eta_\varepsilon \big\|_{L^2} & = \big\| |V_1|^2 - |\Psi(\cdot + a)|^2 
\big\|_{L^2} = \big\| |V_1(\cdot + b - a)|^2 - |\Psi(\cdot + b)|^2 \big\|_{L^2} 
\\
& \leq \big\| |V_1|^2 - |\Psi(\cdot + b)|^2 \big\|_{L^2} + \big\| |V_1(\cdot + b 
- a)|^2 - |V_1|^2 \big\|_{L^2} \\
& \leq d_E(V_1, e^{- i \vartheta} \Psi(\cdot + b)) + C(A) |a - b| \leq (1 + A 
C(A)) d,
\end{align*}
where the continuous function $C(\cdot)$ is provided by 
Lemma~\ref{lem:lipetaV1}. After summation, this yields $\| \varepsilon \|_H + \| 
\eta_\varepsilon \|_{L^2} \leq (1 + A (C(A) + 1)) d$. We require that $(1 + A 
(C(A) + 1)) \rho \leq \delta$, where $\delta$ is given by 
Proposition~\ref{prop:Ecoerc}. This proposition then yields the lower bound
$$
\boE \big( e^{- i \varphi} \Psi(\cdot + a) \big) \geq \kappa \big( \| 
\varepsilon \|_H^2 + \| \eta_\varepsilon \|_{L^2}^2 \big).
$$
To conclude, we finally observe that
$$
\boE \big( e^{- i \varphi} \Psi(\cdot + a) \big) = \boE(\Psi),
$$
by Proposition~\ref{prop:renorm}, and that 
$$
\| \varepsilon \|_H^2 + \| \eta_\varepsilon \|_{L^2}^2 \geq \frac{1}{2} d_E 
\big( V_1, e^{- i \varphi} \Psi(\cdot +a ) \big)^2 \geq \frac{1}{2} d_E \big( 
V_1, \text{Orb}(\Psi) \big)^2,
$$
by definition of the latter. This completes the proof of 
Theorem~\ref{thm:estimcore}, with the choice $\rho = \min \{ 1, \alpha, 
\delta/(1 + A (C(A) + 1)) \}$, and a value of $\kappa$ being half of the 
corresponding value in Proposition~\ref{prop:Ecoerc}.
\end{proof}

\subsection{Concerning orbital stability and Theorem~\ref{thm:stab}}

The proof of Theorem~\ref{thm:stab} assuming Theorem~\ref{thm:estimcore} and 
Proposition~\ref{prop:cauchy} is standard.

\begin{proof}[Proof of Theorem~\ref{thm:stab}]
Let $\Psi_0 \in E$ be such that
$$
d := d_E(V_1, \Psi_0) < \delta,
$$
where $0 < \delta \leq 1$ is a universal constant that will be fixed in the 
course of the proof. For $d = 0$, the conclusion in Theorem~\ref{thm:stab} holds 
since $V_1$ is a stationary solution. Therefore, we assume in the sequel that $d 
> 0$.

First we show that, if $\delta$ is taken smaller than the corresponding value in 
Lemma~\ref{lem:errR}, then we have the estimate
\begin{equation}
\label{eq:boundup}
\boE(\Psi_0) \leq K d^2,
\end{equation}
for some universal constant $K > 0$. Indeed, take $R = \Lambda$, where $\Lambda$ 
is the universal constant provided by the statement of Lemma~\ref{lem:errR}, and 
decompose $\boE(\Psi_0)$ according to~\eqref{eq:decompgood}, i.e.
$$
\boE(\Psi_0) = \frac{1}{2} Q_\Lambda(\varepsilon_0) + N_\Lambda(\varepsilon_0) + 
\frac{1}{2} \boP_\Lambda(\varepsilon_0),
$$
where $\varepsilon_0 := \Psi_0 - V_1.$ By Lemma~\ref{lem:errR}, we obtain
that
$$
|\boP_\Lambda(\varepsilon_0)| \leq K \| \varepsilon_0 \|_H^2 \leq K d^2,
$$
for some universal $K > 0$. A similar estimate holds for 
$Q_\Lambda(\varepsilon_0)$, since $Q_\Lambda$ is a continuous quadratic form on 
$H$ due to Lemma~\ref{lem:comp-emb}. Finally, we check that
$$
|N_\Lambda(\varepsilon_0)| \leq \frac{1}{4} \| \eta_{\varepsilon_0} \|_{L^2}^2 + 
\boI_\Lambda(\varepsilon_0) \leq K \big( \| \eta_{\varepsilon_0} \|_{L^2}^2 + \| 
\varepsilon_0 \|_H^2 \big) = Kd^2,
$$
where $\eta_{\varepsilon_0} := |V_1|^2 - |\Psi_0|^2$, as before.

Define next the constant $C := \max \{ 2, (2 K/\kappa)^{1/2} \}$, where $K > 0$ 
is the constant in~\eqref{eq:boundup}, and $\kappa > 0$ is the constant provided 
by Theorem~\ref{thm:estimcore}. Define then the constant $\delta := \rho/(2 C)$, 
where $\rho > 0$ is also provided by Theorem~\ref{thm:estimcore}.

Let $\Psi_t$ be the solution of the Gross-Pitaevskii equation with initial data 
$\Psi_0$. We pretend that
\begin{equation}
\label{eq:strict}
d_E(V_1, \text{Orb}(\Psi_t)) < C \, d_E(V_1, \Psi_0).
\end{equation}
for any $t \in \R$. Since the map $t \mapsto \Psi_t$ is continuous with values 
into $E$, and since~\eqref{eq:strict} holds for $t = 0$ because $C \geq 2$, it 
suffices to show that the equality
\begin{equation}
\label{eq:absurd}
d_E(V_1, \text{Orb}(\Psi_t)) = C \, d_E(V_1, \Psi_0) 
\end{equation}
for some $t \in \R$ leads to a contradiction. Since $C \, d_E(V_1, \Psi_0) \leq 
C \delta = \frac{\rho}{2}$, we may apply Theorem~\ref{thm:estimcore} to $\Psi_t$ 
when~\eqref{eq:absurd} holds, and conclude that
\begin{equation}
\label{eq:bounddown}
\boE(\Psi_t) \geq \kappa \, d_E(V_1, \text{Orb}(\Psi_t))^2 = \kappa C^2 \, d_E(V_1, 
\Psi_0) ) = 2 K d^2.
\end{equation}
Since $\boE(\Psi_t) = \boE(\Psi_0)$, we deduce from~\eqref{eq:boundup} 
and~\eqref{eq:bounddown} that $K d^2 \geq 2 K d^2$. This is a contradiction 
since $K > 0$ and $d > 0$, and so~\eqref{eq:strict} holds for any $t \in \R$.
\end{proof}

\section{Renormalization of the energy - Proofs of Le\-mma~\ref{lem:errR} and 
Proposition~\ref{prop:renorm} completed}
\label{sect:P1remaining}

We first complete the

\begin{proof}[Proof of Proposition~\ref{prop:renorm}]
We have already shown in the introduction that 
$$
\boE(\Psi) = \lim_{r \to + \infty} \int_{B_r} \big( e_\text{GL}(\Psi) - 
e_\text{GL}(V_1) \big)
$$
is a well-defined quantity when $\Psi \in E$. It remains to prove its invariance 
by translation and phase shift. The latter is immediate, but the former requires 
a short argument.

First, it follows from Lemma~\ref{lem:propprofil} that 
$$
\lim_{r \to + \infty} \int_{B_{R + r} \setminus B_r} e_\text{GL}(V_1) = 0,
$$
for any $R \geq 0$. Let then $\Psi \in E$ and $a \in \R^2$. Since the density 
$e_\text{GL}$ is non-negative, we have
$$
\int_{B_r} e_\text{GL} \big( \Psi(\cdot + a) \big) \leq \int_{B_{r + |a|}} e_\text{GL}(\Psi),
$$
and therefore,
$$
\int_{B_r} \Big( e_\text{GL} \big( \Psi(\cdot + a) \big) - e_\text{GL}(V_1) 
\Big) \leq \int_{B_{r + |a|}} \big( e_\text{GL}(\Psi) - e_\text{GL}(V_1) \big) + 
\int_{B_{r + |a|} \setminus B_r} e_\text{GL}(V_1).
$$
Taking the limit as $r \to + \infty$, we obtain that $\boE(\Psi(\cdot + a)) \leq 
\boE(\Psi)$. It suffices then to interchange the roles of $\Psi$ and $\Psi(\cdot 
+a)$ to obtain the reverse inequality. This completes the proof of 
Proposition~\ref{prop:renorm}.
\end{proof}

We next turn to the end of the

\begin{proof}[Proof of Lemma~\ref{lem:errR}]
Recall that the existence of the quantity 
$$
P_R(\Psi) = \lim_{r \to + \infty} 2 \int_{B_r} (1 - \chi_R)^2 
\frac{x^\perp}{|x|^2} \cdot \langle i \Psi \bar{V}_1, \nabla(\Psi \bar{V}_1) 
\rangle_\C
$$
was already established in the introduction for any function $\Psi \in E$. It 
remains to prove that the existence of universal constants $\delta > 0$, 
$\Lambda > 0$ and $K \geq 1$ such that
\begin{equation}
\label{eq:smallerrR2}
|P_R(\Psi)| \leq \frac{K}{R} \, d_E(\Psi, V_1)^2,
\end{equation}
provided that $d_E(\Psi, V_1) \leq \delta$ and $R \geq \Lambda$. For that 
purpose, we replace the decomposition $\Psi \bar{V}_1 = e^{i\varphi} + w$ 
in~\eqref{eq:PG} by a slight variant, but only available on $B_\Lambda^c$, 
namely
\begin{equation}
\label{eq:PG2}
\Psi \bar{V}_1 = |V_1|^2 \big( e^{i \varphi_\varepsilon} + w_\varepsilon \big),
\end{equation}
where $w_\varepsilon \in H^1(\R^2, \C)$ and $\varphi_\varepsilon \in \dot 
H^1(\R^2, \R)$ will satisfy suitable smallness estimates. Using the property 
that $R \geq \Lambda$, we can modify~\eqref{eq:equiverrR} into
$$
P_R(\Psi) = 2 \int_{\R^2} (1 - \chi_R)^2 |V_1|^4 \frac{x^\perp}{|x|^2} \cdot 
\langle i w_\varepsilon, \nabla w_\varepsilon + 2 i \nabla \varphi_\varepsilon 
e^{i \varphi_\varepsilon} \rangle_\C.
$$
Estimate~\eqref{eq:smallerrR2} then follows from the next 
lemma and the Cauchy-Schwarz inequality, using that
$$
\bigg\| (1 - \chi_R)^2 |V_1|^4 \frac{x^\perp}{|x|^2} \bigg\|_{L^\infty} \leq \frac{K}{R}.
$$
\end{proof}

The decomposition in~\eqref{eq:PG2} is provided by the following lemma.

\begin{lem}
\label{lem:PGalt}
\label{phiw}
There exist constants $K > 0$, $\Lambda > 0$ and $\delta > 0$ such that, given 
any $\Psi = V_1 + \varepsilon \in E$ with
$$
\| \eta_\varepsilon \|_{L^2} + \| \nabla (\varepsilon \bar{V}_1) \|_{L^2} \leq 
\delta,
$$
there exist $\varphi_\varepsilon \in \dot{H}^1(\R^2, \R)$ and $w_\varepsilon \in 
H^1 (\R^2, \C)$ such that
$$
\Psi \bar{V}_1 = |V_1|^2 \big( e^{i \varphi_\varepsilon} + w_\varepsilon \big)
$$
on $B_\Lambda^c$. Moreover, the functions $\varphi_\varepsilon$ and 
$w_\varepsilon$ satisfy
$$
\| w_\varepsilon \|_{H^1(\R^2)} \leq K \big(\| \nabla (\varepsilon \bar{V}_1) 
\|_{L^2} + \| \eta_\varepsilon \|_{L^2} \big),
$$
and
$$
\| \nabla \varphi_\varepsilon \|_{L^2 (\R^2)} \leq K \| \nabla (\varepsilon 
\bar{V}_1) \|_{L^2}.
$$
\end{lem}

\begin{proof}
In the whole proof, the notation $A \lesssim B$, for any arbitrary non-negative 
numbers $A$ and $B$, means that $A \leq K B$ for some universal constant $K$, 
whose exact value is irrelevant for the argument. We denote by $\chi(D)$ the 
cut-off in Fourier space induced by the cut-off function $\chi$. Recall that 
this function is radial, with $\chi \equiv 1$ in $B_1$ and $\chi \equiv 0$ 
outside $B_2$, so that $\chi(D)$ only keeps the small frequencies. For $\Lambda 
> 1$ and $\Psi = V_1 + \varepsilon$, we decompose the function $\varepsilon 
\bar{V}_ 1$ on $B_\Lambda^c$ as
$$
\varepsilon \bar{V}_1 = \chi (D) \big( \varepsilon \bar{V}_1 \big) + \frac{(1 - 
\chi (D)) \big( \varepsilon \bar{V}_1 \big)}{|V_1|^2} |V_1|^2.
$$
We let $w_0 := (1 - \chi (D)) \big( \varepsilon \bar{V}_1 \big)/|V_1|^2$ on 
$B_\Lambda^c$, and check using Lemma~\ref{lem:propprofil} that
$$
\| w_0 \|_{H^1(B_\Lambda^c)} \lesssim \| \nabla (\varepsilon \bar{V}_1) 
\|_{L^2}.
$$
We now aim to prove that, provided that $\Lambda > 1$ is large enough and $0 < 
\delta \leq 1$ is small enough, the function $w := \chi (D) \big( \varepsilon 
\bar{V}_1 \big) + |V_1|^2$ satisfies $|w| \geq 1/2$ on $B_\Lambda^c$. For that 
purpose, we decompose
\begin{equation}
\label{Estw}
1 - |w|^2 = 1 - \big| |V_1|^2 + \varepsilon \bar{V}_1 \big|^2 \, - \, \big| (1 
- \chi (D)) \big( \varepsilon \bar{V}_1 \big) \big|^2 \, + \, 2 \big\langle (1 - 
\chi (D)) \big( \varepsilon \bar{V}_1 \big), |V_1|^2 + \varepsilon \bar{V}_1 
\big\rangle_\C,
\end{equation}
and we start by estimating these three terms in $L^2(B_\Lambda^c) $. For the 
first one, we rewrite
$$
1 - \big| |V_1|^2 + \varepsilon \bar{V}_1 \big|^2 = 1 - |V_1|^4 - 2 \langle 
\varepsilon, V_1 \rangle_\C |V_1|^2 - |V_1|^2 |\varepsilon|^2 = 1 - | V_1 |^4 + 
\eta_\varepsilon |V_1|^2,
$$
and therefore,
$$
\big\| 1 - \big| |V_1|^2 + \varepsilon \bar{V}_1 \big|^2 
\big\|_{L^2(B_\Lambda^c)} \leq \| 1 - |V_1|^4 \|_{L^2(B_\Lambda^c)} + \| 
\eta_\varepsilon \|_{L^2 (B_\Lambda^c)} \lesssim \Lambda^{- 1} + \delta.
$$
Here, we have used Lemma~\ref{lem:propprofil} in order to estimate the decay 
with respect to $\Lambda$. For the second term, we infer from the Sobolev 
embedding theorem that
$$
\Big\| \big| (1 - \chi (D)) \big( \varepsilon \bar{V}_1 \big) \big|^2 
\Big\|_{L^2} \lesssim \big\| (1 - \chi (D)) \big( \varepsilon \bar{V}_1 
\big) \big\|_{H^1}^2 \lesssim \| \nabla (\varepsilon \bar{V}_1) \|_{L^2 
}^2 \leq \delta^2 \leq \delta.
$$
Finally, since $\big| |V_1|^2 + \varepsilon \bar{V}_1 \big| \lesssim 1 + |\Psi| 
\lesssim 1 + |\eta_\varepsilon|^{1/2}$ and $|\eta_\varepsilon|^{1/2} \in 
L^4(\R^2)$, we deduce again from the Sobolev embedding theorem that the third 
term is controlled by
\begin{align*}
& \Big\| \big\langle (1 - \chi (D)) \big( \varepsilon \bar{V}_1 \big), |V_1|^2 + 
\varepsilon \bar{V}_1 \big\rangle_\C \Big\|_{L^2}\\
\lesssim & \big\| (1 - \chi (D)) \big( \varepsilon \bar{V}_1 \big) 
\big\|_{L^2} + \big\| (1 - \chi (D)) \big( \varepsilon \bar{V}_1 \big) 
\big\|_{H^1} \big\| \eta_\varepsilon \big\|_{L^2}^\frac{1}{2}\\
\lesssim & \| \nabla (\varepsilon \bar{V}_1) \|_{L^2} \big( 1 + \| 
\eta_\varepsilon \|_{L^2}^\frac{1}{2} \big) \lesssim \delta \big( 1 + 
\delta^\frac{1}{2} \big) \lesssim \delta.
\end{align*}
Combining these estimates in~\eqref{Estw}, we deduce that
$$
\| 1 - |w|^2 \|_{L^2(B_\Lambda^c)} \lesssim \Lambda^{- 1} + \delta.
$$
 
To obtain a uniform control rather than an $L^2$-one, we rely on the 
Gagliardo-Nirenberg type inequality
$$
\big\| 1 - |w|^2 \big\|_{L^\infty(B_\Lambda^c)} \lesssim \big\| 1 - | w|^2 
\big\|_{L^2(B_\Lambda^c)}^\frac{1}{2} \big\|D^2 (| w|^2) 
\big\|_{L^2(B_\Lambda^c)}^\frac{1}{2},
$$
where the proportionality constant does not depend on $\Lambda \geq 1$, as it 
can be seen by scaling from $B_1^c$ and invoking a standard Sobolev extension 
theorem. Then we bound
$$
\big\| D^2 (| w|^2) \big\|_{L^2(B_\Lambda^c)} \lesssim \big\|D^2 w 
\big\|_{L^\infty(B_\Lambda^c)} \big\| 1- |w|^2 \big\|_{L^2(B_\Lambda^c)} +
\big\| D^2 w \big\|_{L^2(B_\Lambda^c)} + \big\| D w \big\|_{L^4(B_\Lambda^c)}^2.
$$
From Lemma~\ref{lem:propprofil} and the fact that $\chi(D)(\eps \bar V_1)$ only 
has small frequencies, we also have
\begin{equation}
\label{eq:pentecote}
\| D^k w \|_{L^2(B_\Lambda^c)} \lesssim \Lambda^{-1} + \delta,
\end{equation}
with a constant depending only on $k \geq 1$. Combining the previous inequalities and
the Sobolev embedding theorem, we deduce
$$
\big\| 1 - |w|^2 \big\|_{L^\infty(B_\Lambda^c)} \lesssim \Lambda^{-1} + \delta.
$$
Therefore, we obtain $|w| \geq 1/2$ on $B_\Lambda^c$ as claimed, provided that 
$\Lambda > 0$ is chosen sufficiently large and $\delta > 0$ sufficiently small. 
Under these conditions, there exists a function $\varphi_\varepsilon \in \dot 
H^1(B_\Lambda^c)$ such that
$$
\chi (D) \big( \varepsilon \bar{V}_1 \big) + | V_1 |^2 = w = \big| \chi (D) 
\big( \varepsilon \bar{V}_1 \big) + |V_1|^2 \big| \, e^{i \varphi_\varepsilon}.
$$
We decompose $\chi (D) \big( \varepsilon \bar{V}_1 \big) = |V_1|^2 \big( e^{i 
\varphi_\varepsilon} - 1 \big) + |V_1|^2 w_1$, with
$$
w_1 := \frac{1}{| V_1 |^2} (| | V_1 |^2 + \chi (D) (\varepsilon\bar{V}_1) | - | 
V_1 |^2) e^{i \varphi_\varepsilon} = \frac{1}{| V_1 |^2} (| w | - | V_1 |^2) 
e^{i \varphi_\varepsilon}.
$$
Since $|w| \geq 1/2$ on $B_\Lambda^c$, we observe that $\big| |w| - |V_1|^2 
\big| \lesssim \big| |w|^2 - |V_1|^4 \big|$ on $B_\Lambda^c$, which implies that
$$
\big\| w_1 \big\|_{L^2(B_\Lambda^c)} \lesssim \big\| |w|^2 - |V_1|^4 \big\|_{L^2 
(B_\Lambda^c)} \lesssim \| \eta_\varepsilon \|_{L^2} + \| \nabla 
(\varepsilon\bar{V}_1) \|_{L^2}.
$$
We next write
$$
\| \nabla w_1 \|_{L^2 (B_\Lambda^c)} \lesssim \| \nabla (w_1 e^{- i 
\varphi_\varepsilon}) \|_{L^2 (B_\Lambda^c)} + \| w_1 \|_{L^\infty(B_\Lambda^c)} 
\| \nabla \varphi_\varepsilon \|_{L^2 (B_\Lambda^c)}.
$$
We first observe that
$$
\nabla \big( w_1 e^{- i \varphi_\varepsilon} \big) = |V_1|^2 \nabla \Big( 
\frac{1}{|V_1|^2} \Big) w_1 e^{- i \varphi_\eps} + \frac{1}{|V_1|^2} \nabla 
\chi(D) (\eps\bar V_1).
$$
For $\Lambda > 1$, we therefore obtain
\begin{align*}
\big\| \nabla (w_1 e^{- i \varphi_\varepsilon}) \big\|_{L^2(B_\Lambda^c)} & 
\lesssim \| w_1 \|_{L^2(B_\Lambda^c)} + \| \nabla \chi (D) \big( \varepsilon 
\bar{V}_1) \big\|_{L^2(B_\Lambda^c)}\\
& \lesssim \| \nabla (\varepsilon\bar{V}_1)\|_{L^2} + \| 
\eta_\varepsilon\|_{L^2}.
\end{align*}
We also have
$$
\| w_1 \|_{L^\infty(B_\Lambda^c)} \lesssim \| 1 - |w|^2 
\|_{L^\infty(B_\Lambda^c)} + 1 \lesssim 1.
$$
Finally, since $e^{i \varphi_\varepsilon} = w/|w|$, with $|w| \geq 1/2$ on 
$B_\Lambda^c$, we have
\begin{align*}
\| \nabla \varphi_\varepsilon \|_{L^2(B_\Lambda^c)} & \lesssim \big\| \langle i 
w, \nabla w \rangle_\C \big\|_{L^2(B_\Lambda^c)} \\
& = \Big\| \big\langle i \big( \chi (D) (\varepsilon \bar{V}_1) + |V_1|^2 
\big), \nabla \big( \chi (D) (\varepsilon \bar{V}_1) + |V_1|^2 \big) 
\big\rangle_\C \Big\|_{L^2(B_\Lambda^c)}
\end{align*}
and the latter expression is controlled by 
$$
\big\| \chi (D) (\varepsilon \bar{V}_1) \big\|_{L^2(B_\Lambda^c)} \Big( 1 + 
\big\| \chi (D) (\varepsilon \bar{V}_1) \big\|_{L^\infty(B_\Lambda^c)} \Big) 
\lesssim \| \nabla (\varepsilon \bar{V}_1) \|_{L^2},
$$
since $\delta \leq 1$. This leads to
$$
\| \nabla w_1 \|_{L^2 (B_\Lambda^c)} \lesssim \| \nabla (\varepsilon \bar{V}_1) 
\|_{L^2} + \| \eta_\varepsilon \|_{L^2}.
$$
We define $w_\varepsilon := w_0 + w_1$, which satisfies
$$
\| w_\varepsilon \|_{H^1 (B_\Lambda^c)} \lesssim \| \nabla (\varepsilon \bar{V}_1) 
\|_{L^2} + \| \eta_\varepsilon \|_{L^2}.
$$
We check that, by construction, $\Psi \bar{V}_1 = |V_1|^2 \big( e^{i 
\varphi_\varepsilon} + w_\varepsilon \big)$ on $B_\Lambda^c$. It suffices then 
to extend $w_\eps$ and $\varphi_\eps$ to $\R^2$ using the Sobolev extension 
theorem in order to complete the proof.
\end{proof}

\section{Minimality - Proof of Proposition~\ref{prop:minim} completed} 
\label{sect:P2remaining}

In the introduction, we have already mentioned that in order to prove 
Proposition~\ref{prop:minim}, it suffices to show that the renormalized 
Ginzburg-Landau energy is non-negative. For that purpose, we invoke results 
from~\cite{Mirones1} regarding the asymptotics of the Ginzburg-Landau energy for 
functions on a fixed bounded domain with fixed boundary data.

Let $(R_n)_{n \in \N}$ be a sequence of positive numbers such that $R_n \to + 
\infty$. Define $\varepsilon_n := R_n^{- 1}$ and let $u_n(x) := \Psi(x / 
\varepsilon_n)$ on $B_1$. By scaling, we have
$$
\int_{B_{R_n}} e_\text{GL}(\Psi) = \int_{B_1} e_{GL_{\varepsilon_n}}(u_n) := 
\int_{B_1} \Big( \frac{|\nabla u_n|^2}{2} + \frac{(1 - |u_n|^2 )^2}{4 
\varepsilon_n^2} \Big),
$$
and therefore we only need to show that
$$
\liminf_{n \to + \infty} \int_{B_1} \bigg( e_{GL_{\varepsilon_n}}(u_n) - 
e_{GL_{\varepsilon_n}} \Big( V_1 \Big( \frac{\cdot}{\varepsilon_n} \Big) \Big) 
\bigg) \geq 0.
$$
Since the boundary data of $u_n$ on $\partial B_1$ is not fixed, we cannot 
invoke the results of~\cite{Mirones1} directly. Instead, we choose the sequence 
$(R_n)_{n \in \N}$ so that $u_n$ has almost minimal energy on $\partial B_1$, 
and then perform an extension to a slightly larger domain with a fixed boundary 
data. More precisely, since $\nabla (\Psi \bar{V}_1)$ and $1 - |\Psi|^2$ both 
belong to $L^2(\R^2)$, we may find a sequence $(R_n)_{n \in \N}$ such that $R_n 
\to + \infty$ and
\begin{equation}
\label{eq:petit}
\int_{\partial B_{R_n}} \Big( |\nabla(\Psi \bar{V}_1)|^2 + (1 - |\Psi|^2)^2 
\Big) = o \Big( \frac{1}{R_n} \Big),
\end{equation}
as $n \to + \infty$. For $n$ sufficiently large, it follows that we may write 
$$
u_n \big( e^{i \theta} \big) = \big| u_n(e^{i \theta}) \big| \, e^{i \big( 
\theta + \phi_n(\theta) \big)} = \rho_n(\theta) \, e^{i \big( \theta + 
\phi_n(\theta) \big)}
$$
on $\partial B_1$, where
\begin{equation}
\label{eq:petit2}
\int_0^{2 \pi} \Big( \frac{|\partial_\theta \rho_n|^2}{2} + \rho_n^2 
\frac{|\partial_\theta \phi_n|^2}{2} + \frac{(1 - \rho_n^2)^2}{4 
\varepsilon_n^2}\Big) \, d\theta = o(1),
\end{equation}
as $n \to + \infty$. Let $\bar{\phi}_n$ denote the mean of $\phi_n$ on
$\partial B_1$. Fix $\delta > 0$ and consider the extension 
$$
u_n \big( r e^{i \theta} \big) = \Big( \lambda(r) \rho_n(\theta) +
1 - \lambda(r) \Big) e^{i \big( \theta + \lambda(r) \phi_n(\theta) + (1 -
\lambda(r)) \bar{\phi}_n \big)}
$$
of $u_n$ for $1 \leq r \leq 1+ \delta$, where $\lambda(r) := (1 + \delta - 
r)/\delta$. We deduce from~\eqref{eq:petit2} and elementary computations that 
the function $v_n := e^{- i \bar{\phi}_n} u_n$, which is defined on $B_{1 + 
\delta}$ and with fixed boundary data $v_n \big( (1 + \delta) e^{i \theta} \big) 
= e^{i \theta}$ on $\partial B_{1 + \delta}$, satisfies
$$
\int_{B_{1 + \delta}} e_{GL_{\varepsilon_n}}(v_n) \leq \int_{B_1} 
e_{GL_{\varepsilon_n}}(u_n) + \pi \log(1 + \delta) + o(1),
$$
as $n \to + \infty$. By~\cite[Corollaire 2]{Mirones1} and the remark following 
it, we have
$$
\liminf_{n \to + \infty} \int_{B_{1 + \delta}} \bigg( 
e_{GL_{\varepsilon_n}}(v_n) - e_{GL_{\varepsilon_n}} \Big( V_1 \Big( 
\frac{\cdot}{\varepsilon_n} \Big) \Big) \bigg) \geq 0,
$$
and therefore
$$
\liminf_{n \to + \infty} \int_{B_1} \bigg( e_{GL_{\varepsilon_n}}(u_n) - 
e_{GL_{\varepsilon_n}} \Big( V_1 \Big( \frac{\cdot}{\varepsilon_n} \Big) \Big) 
\bigg) \geq - \pi \log(1+\delta).
$$
The conclusion follows letting $\delta \to 0$. \qed

\section{Second order expansion - Proof of Lemma~\ref{lem:decomp}}

First recall the decomposition in~\eqref{eq:decompboE} of the energy
\begin{align*}
\boE(\Psi) = & \frac{1}{2} \| \Psi \|_H^2 - \int_{\R^2} |\nabla V_1|^2 (1 - 
|\Psi|^2)\\
& - \int_{\R^2} \Big\langle \nabla(\Psi \bar{V}_1), \Big( \nabla 
\bar{V}_1 + i (1 - \chi_R)^2 \frac{x^\perp}{|x|^2} \bar{V}_1 \Big) \Psi 
\Big\rangle_\C + \frac{1}{2} P_R(\Psi)\\
& + \int_{\R^2} \frac{1}{4} \big( (1 - |\Psi|^2)^2 - (1 - |V_1|^2)^2 \big).
\end{align*}
We write $\Psi = V_1 + \varepsilon \in E$ and then expand the various terms in the 
previous decomposition. Note that no use is made of the equation~\eqref{eq:GL} 
satisfied by $V_1$ until the very end of the argument. The qualitative 
properties of $V_1$ in Lemma~\ref{lem:propprofil} associated to the embedding 
given by Lemma~\ref{lem:comp-emb} for $\varepsilon \in H$ are sufficient to justify all 
computations.

After elementary algebra and integration by parts, we obtain
$$
\frac{1}{2} \| \Psi \|_H^2 = \frac{1}{2} \| \varepsilon \|_H^2 + \int_{\R^2} \big( 
T_{0a} + \langle T_{1a}, \varepsilon \rangle_\C \big),
$$
where 
$$
T_{0a} := \frac{1}{2} \big| \nabla (|V_1|^2) \big|^2 + \frac{1}{2} |\nabla 
V_1|^2 (1 - | V_1 |^2),
$$
and 
$$
T_{1a} := - \Delta \big( | V_1 |^2 \big) V_1 - \text{div} \big( (1 - | V_1 |^2) 
\nabla V_1 \big).
$$
Using that $|V_1|$ and $\chi_R$ are radial functions, we also obtain
\begin{align*}
& - \int_{\R^2} \Big\langle \nabla (\Psi \bar{V}_1), \Big( \nabla \bar{V}_1 + i 
(1 - \chi_R)^2 \frac{x^\perp}{|x|^2} \bar{V}_1 \Big) \Psi \Big\rangle_\C \\
= & - \int_{\R^2} \Big\langle \nabla (\varepsilon \bar{V}_1), \Big( \nabla 
\bar{V}_1 + i (1 - \chi_R)^2 \frac{x^\perp}{|x|^2} \bar{V}_1 \Big) \varepsilon 
\Big\rangle_\C + \int_{\R^2} \big( T_{0b} + \langle T_{1b}, \varepsilon 
\rangle_\C \big),
\end{align*}
where
$$
T_{0b} := - \frac{1}{2} \big| \nabla(|V_1|^2) \big|^2 \quad \text{and} \quad 
T_{1b} := \Delta \big( |V_1|^2 \big) V_1 - \text{div} \big( |V_1|^2 \nabla V_1 
\big) - V_1 |\nabla V_1|^2.
$$
Finally, recalling the definition $\eta_\varepsilon = - 2 \langle \varepsilon, 
V_1 \rangle_\C - |\varepsilon|^2$ in~\eqref{eq:defetaeps}, we have
$$
- \frac{1}{2} |\nabla V_1|^2 (1 - |\psi|^2) = T_{0c} + \langle T_{1c}, 
\varepsilon \rangle_\C + \frac{1}{2} |\nabla V_1|^2 |\varepsilon|^2,
$$
where
$$
T_{0c} := - \frac{1}{2} |\nabla V_1|^2 (1 - |V_1|^2) \quad \text{and} \quad 
T_{1c} := |\nabla V_1|^2 V_1.
$$
Similarly, we compute
$$
\frac{1}{4} \big( (1 - |\psi|^2)^2 - (1 - |V_1|^2)^2 \big) = \langle T_{1d}, 
\varepsilon \rangle_\C - \frac{1}{2} (1 - |V_1|^2) |\varepsilon|^2 + \frac{1}{4} 
\eta_{\varepsilon}^2,
$$
where
$$
T_{1d} := - (1 - |V_1 |^2) V_1.
$$
It is immediate that $T_{0a} + T_{0b} + T_{0c} = 0$, and we check that
$$
T_{1a} + T_{1b} + T_{1c} + T_{1d} = - \Delta V_1 - (1 - |V_1|^2) V_1 = 0,
$$
since $V_1$ is a solution to~\eqref{eq:GL}. Lemma~\ref{lem:decomp} then follows 
after summing the previous identities. \qed

\section{Quadratic form - Proof of Proposition~\ref{prop:quad}}
\label{sect:quad}

In this section, we establish the coercivity properties of the quadratic form 
$Q_R$ as stated in Proposition~\ref{prop:quad}. In particular, we assume 
throughout this section that $R \ge 1$. For $\varepsilon \in H$, recall that
$$
Q_R(\varepsilon) := \boQ_R(\varepsilon) + 2\boI_R(\varepsilon)
$$
where
\begin{align*}
\boQ_R(\varepsilon) = & \big\| \varepsilon \big\|_H^2 - \int_{\R^2} \big( 1 - 
|V_1|^2 - |\nabla V_1|^2 \big) \, |\varepsilon|^2\\ & - 2 \int_{\R^2} 
\Big\langle \nabla(\varepsilon \bar{V}_1), \Big( \nabla \bar{V}_1 + i 
\frac{x^\perp}{|x|^2} (1 -\chi_R)^2 \bar{V}_1 \Big) \varepsilon \Big\rangle_\C,
\end{align*}
and
$$
\boI_R(\varepsilon) = \int_{\R^2} \chi_R^2 \big\langle \varepsilon, V_1 
\big\rangle_\C^2.
$$

Combining the decay and regularity properties of $V_1$ stated in 
Lemma~\ref{lem:propprofil} with the embeddings for the Hilbert space $H$ in 
Lemma~\ref{lem:comp-emb}, we can check that the quadratic form $Q_R$ is 
well-defined and continuous on $H$. Moreover, its associated self-adjoint 
operator is a compact perturbation of the identity on $H$.

In the course of the proof, we shall use multiple times the following elementary 
consequence of the Fredholm theory.

\begin{lem}
\label{lem:fred}
Let $(X, \| \cdot \|_X)$ be a real Hilbert space and $Q$ be a real continuous 
quadratic form on $X$. Assume that, for an inner product on $X$ whose norm is 
equivalent to the original one, $Q$ is associated to a compact perturbation of 
the identity on $X$. Assume also that $Q(x) > 0$ for all $x \in X \setminus \{ 0 
\}$. Then there exists $\kappa_0 > 0$ such that
$$
Q(x) \geq \kappa_0 \| x \|_X^2,
$$
for any $x \in X$.
\end{lem}

The Hilbert space $H$ is naturally decomposed into orthogonal Fourier sectors 
through the formula~\footnote{The shift in the exponent is convenient for later 
symmetry purposes.}
$$
\varepsilon(r, \theta) = \sum_{j \in \Z} \varepsilon_j(r) e^{i (j + 1) \theta}.
$$
Elementary computations show that
$$
\big\| \varepsilon_j e^{i (j + 1) \theta} \big\|_H \simeq \| \varepsilon_j \|_{H_j},
$$
with universal similarity constants independent of $j \in \Z$. Here, the Hilbert 
space $H_j$ is defined through the norm
\begin{equation}
\label{eq:defhj}
\big\| e \big\|_{H_j}^2 := \int_0^{+ \infty} \bigg( |e'(r)|^2 + \Big( 
\frac{j^2}{1 + r^2} + \frac{(j + 1)^2}{r^2 (1 + r^2)} \Big) |e(r)|^2 \bigg) r \, 
dr.
\end{equation}

In this framework, the quadratic forms $\boQ_R$ and $\boI_R$ may be decomposed 
as
\begin{equation}
\label{eq:decompboQ}
\boQ_R(\varepsilon) = \sum_{j \in \Z} \Big( Q_{R, j}(a_j) + Q_{R, j}(b_j) \Big),
\end{equation}
and
\begin{equation}
\label{eq:decompboI}
\boI_R(\varepsilon) = I_R(a_0) + \frac{1}{2} \sum_{j = 1}^{+ \infty} \Big( I_R 
\big( a_j + a_{- j} \big) + I_R \big( b_j - b_{- j} \big) \Big),
\end{equation}
where we have set $\varepsilon_j =: a_j + i b_j$, with $a_j$ and $b_j$ being 
real-valued functions. In the previous identities, the quadratic forms $Q_{R, 
j}$ and $I_R$ are given by
$$
Q_{R, j}(e) := \int_0^{+ \infty} \bigg( |e'|^2 + \Big( \frac{(j + 1)^2}{r^2} -
2 j \frac{(1 - \chi_R)^2}{r^2} \rho_1^2 - (1 - \rho_1^2) \Big) |e|^2 \bigg) r \, dr,
$$
and
$$
I_R(e) := \int_0^{+ \infty} \rho_1^2 \chi_R^2 |e|^2 \, r \, dr.
$$
We now describe the coercivity properties of the quadratic forms $Q_{R, j}$ and 
$I_R$, as well as of suitable combinations of them, beginning by $Q_{R, 0}$.

\begin{lem}
\label{lem:Q0pos}
The quadratic form $Q_{R, 0}$ is independent of the number $R$ and it satisfies
\begin{equation}
\label{eq:thebes1}
Q_0(e) := Q_{R, 0}(e) = \int_0^{+ \infty} \rho_1^2 \Big| \Big( \frac{e}{\rho_1} 
\Big)' \Big|^2 \, r \, dr \geq 0.
\end{equation}
for any $e \in H_0$.
\end{lem}

\begin{proof}
By definition, the quadratic form $Q_{R, 0}$ does not depend on the number $R$. 
Moreover, it follows from the Leibniz rule that
$$
\int_0^{+ \infty} \rho_1^2 \Big| \Big( \frac{e}{\rho_1} \Big)' \Big|^2 \, r \, 
dr = \int_0^{+ \infty} \Big( |e'|^2 + \frac{(\rho_1')^2}{\rho_1^2} |e|^2 - 
(|e|^2)' \frac{\rho_1'}{\rho_1} \Big) r \, dr.
$$
Integrating by parts, we obtain
$$
- \int_0^{+ \infty} (|e|^2)' \frac{\rho_1'}{\rho_1} \, r \, dr = \int_0^{+ \infty}
\Big( \frac{\rho_1''}{\rho_1} + \frac{\rho_1'}{r \rho_1} \Big) |e|^2\, r \, dr.
$$
Indeed, the boundary terms vanish. This follows from Lemma~\ref{lem:propprofil}, 
as well as the fact that a function $e \in H_0$ satisfies $e(0) = 0$ and has an 
at most logarithmic growth at infinity. In particular, this latter property 
guarantees that $|e(r)|^2/r^2 = o(1)$ as $r \to + \infty$. 
Identity~\eqref{eq:thebes1} then follows from using the equation~\eqref{eq:rho1} 
for the function $\rho_1$.
\end{proof}

We next compare the quadratic forms $Q_{R, j}$ with $Q_0$.

\begin{lem}
\label{lem:QjQ0}
Let $j \in \Z \setminus \{ - 1 \}$ and $e \in H_j \subset H_0$. We have
$$
Q_{R, j}(e) \geq Q_0(e) \geq 0.
$$
More precisely, for $j \neq - 2$, we have
\begin{equation}
\label{eq:thebes3}
Q_{R, j}(e) - Q_0(e) \geq \frac{1}{3} \int_0^{+ \infty} \frac{j^2}{r^2} |e|^2 \, 
r \, dr,
\end{equation}
while for $j = - 2$,
\begin{equation}
\label{eq:thebes4}
Q_{R, - 2}(e) - Q_0(e) = 4 \int_0^{+ \infty} \rho_1^2 (1 - \chi_R)^2 
\frac{|e|^2}{r^2} \, r \, dr.
\end{equation}
\end{lem}

\begin{proof}
We compute
$$
Q_{R, j}(e) - Q_0(e) = \int_0^{+ \infty} \Big( j^2 + 2 j \big( 1 -
(1 - \chi_R)^2 \big) \rho_1^2 \Big) \frac{|e|^2}{r^2} \, r \, dr,
$$
which gives~\eqref{eq:thebes4} for $j = - 2$. For $j \geq 0$, the conclusion 
follows from the inequality $j^2 + 2 j (1 - (1 - \chi_R)^2) \rho_1^2 \geq j^2$. 
For $j \leq - 3$, we instead write $j^2 + 2 j (1 - (1 -\chi_R)^2) \rho_1^2 \geq 
j^2 + 2 j$, and then $j^2 + 2 j \geq j^2/3$.
\end{proof}

As a consequence of the formula~\eqref{eq:thebes1} for $Q_0$, we obtain the 
following coercivity result for this quadratic form.

\begin{cor}
\label{cor:niveau0}
There exists $\kappa_0 > 0$ such that
$$
Q_0(e) + 2I_R(e) \geq \kappa_0 \Big( \| e \|_{H_0}^2 + I_R(e) \Big),
$$
for any $e\in H_0$. Moreover, if $e$ satisfies the orthogonality condition
$$
\int_0^{+ \infty} \chi(r) e(r) \rho_1(r) \, r \, dr = 0,
$$
then we have 
$$
Q_0(e) \geq \kappa_0 \| e \|_{H_0}^2.
$$
\end{cor}

\begin{proof}
We deduce from Lemma~\ref{lem:Q0pos} that the quadratic form $e \mapsto Q_0(e) + 
I_1(e)$ is positive definite. Moreover, it is as $Q_0$ associated to a compact 
perturbation of the identity in $H_0$. Therefore, there exists some constant 
$\kappa > 0$, independent of $R$, such that
$$
Q_0(e) + I_R(e) \geq Q_0(e) + I_1(e) \geq \kappa \| e \|_{H_0}^2,
$$
so that
$$
Q_0(e) + 2 I_R(e) \geq \kappa \| e \|_{H_0}^2 + I_R(e). 
$$
Similarly, the quadratic form $e \mapsto Q_0(e) + (\int_0^{+ \infty} \chi(r) 
e(r) \rho_1(r) \, r \, dr)^2$ is positive definite and associated to a compact 
perturbation of the identity. The conclusion follows for some suitable $\kappa_0 
\leq \min \{ \kappa, 1 \}$.
\end{proof}

Similarly, we derive from~\eqref{eq:thebes1} and~\eqref{eq:thebes3} the 
following coercivity result for the quadratic forms $Q_{R, j}$, with $j \neq - 
2, - 1, 0$.

\begin{cor}
\label{cor:jlarger3}
There exists $\kappa_0 > 0$, independent of $R$, such that
$$
Q_{R, j}(e) \geq \kappa_0 \| e \|_{H_j}^2,
$$
for any $j\in \Z \setminus \{ - 2, - 1, 0 \}$ and any $e \in H_j$.
\end{cor}

\begin{proof}
As for the proof of Corollary~\ref{cor:niveau0}, the quadratic form $e \mapsto 
Q_0(e) + \frac{1}{4} \int_0^1\rho_1^2 |e|^2 \, r \, dr$ is positive definite. Moreover, 
it is as $Q_0$ associated to a compact perturbation of the identity in $H_0$. 
Therefore, we infer the existence of $\kappa_0 >0$ such that
\begin{equation}
\label{eq:thebes2}
Q_0(e) \geq \kappa_0 \big\| e \big\|_{H_0}^2 - \frac{1}{4} \int_0^1 \rho_1^2 
|e|^2 \, r \, dr,
\end{equation}
for any $e \in H_0$. Since $H_j \subset H_0$ for $j \neq - 1$, the summation 
of~\eqref{eq:thebes3} and~\eqref{eq:thebes2} yields the conclusion for $j 
\neq - 2$, and for some possibly smaller value of $\kappa_0$, but that can be 
chosen independently of $j$.
\end{proof}

We next consider the more tedious case $j = - 2$.

\begin{cor}
\label{cor:niveau2}
There exists $\kappa_0 > 0$, independent of $R$, such that we have
$$
Q_{R, 2}(e) + Q_{R, - 2}(f) + I_R(e \pm f) \geq \kappa_0 \Big( \| e \|_{H_2}^2 + 
\| f \|_{H_{- 2}}^2 + I_R(e \pm f) \Big),
$$
for any $e \in H_2$ and $f \in H_{- 2}$.
\end{cor}

\begin{proof}
We derive from Lemma~\ref{lem:QjQ0} that
\begin{equation}
\label{eq:ici0}
Q_{R, 2}(e) + Q_{R, - 2}(f) \geq Q_0(e) + Q_0(f) + \int_0^{+ \infty} 
\frac{|e|^2}{r^2} \, r \, dr + \int_0^{+ \infty} \rho_1^2 (1 -\chi_R)^2 
\frac{|f|^2}{r^2} \, r \, dr.
\end{equation}
Since $0 \leq \rho_1 \leq 1$, we can split
$$
\int_0^{+ \infty} \frac{|e|^2}{r^2} \, r \, dr \geq \int_0^{+ \infty} \rho_1^2 
\chi_R^2 \frac{|e|^2}{r^2} \, r \, dr + \int_0^{+ \infty} \rho_1^2 (1 - 
\chi_R)^2 \frac{|e|^2}{r^2} \, r \, dr.
$$
After summation of~\eqref{eq:ici0} with the inequality
$$
I_R(e \pm f) = \int_0^{+ \infty} \rho_1^2 \chi_R^2 |e \pm f|^2 \, r \, dr \geq
\int_1^\infty \rho_1^2 \chi_R^2 \frac{|e \pm f|^2}{r^2} \, r \, dr,
$$
we obtain the estimate 
\begin{equation}
\label{eq:thebes5}
Q_{R, 2}(e) + Q_{R, - 2}(f) + I_R(e \pm f) \geq Q_0(e) + Q_0(f) + \frac{1}{6} 
\int_1^{+ \infty} \rho_1^2 \Big( \frac{|e|^2}{r^2} + \frac{|f|^2}{r^2}\Big) \, r 
\, dr,
\end{equation}
using that $e^2 + (e \pm f)^2 \geq e^2/2 + f^2/3$ and $\chi_R^2 + (1 - \chi_R)^2 
\geq 1/2$. The quadratic form in the right-hand side of~\eqref{eq:thebes5} is 
positive definite by Lemma~\ref{lem:Q0pos} and independent of $R$. It is associated to 
compact perturbations of the identity for norms that are equivalent to the ones 
in $H_2$ and $H_{- 2}.$ It follows that
$$
Q_{R, 2}(e) + Q_{R, - 2}(f) + I_R(e \pm f) \geq \kappa \Big( \| e \|_{H_2}^2 + 
\| f \|_{H_{- 2}}^2 \Big),
$$
for some $\kappa > 0$. We conclude by observing that
$$
Q_{R, 2}(e) + Q_{R, - 2}(f) + I_R(e \pm f) \geq \frac{1}{2} \Big( Q_{R, 2}(e) + 
Q_{R, - 2}(f) + I_R(e \pm f) \Big) + \frac{1}{2} I_R(e \pm f),
$$
by Lemma~\ref{lem:QjQ0}, and therefore choosing $\kappa_0 = \min \{ \kappa, 1 \}/2$.
\end{proof}

We finally establish some coercivity for the quadratic forms $Q_{R, \pm 1}$ 
under suitable orthogonality conditions.

\begin{prop}
\label{prop:H1Hm1}
There exist $\kappa_0 > 0$, $R_0 \geq 2$ and $C > 0$ such that, given any $R 
\geq R_0$ and any real-valued functions $e \in H_1$ and $f \in H_{- 1}$ 
satisfying the orthogonality conditions
\begin{equation}
\label{eq:ortho1m1}
\int_0^{+ \infty} \chi \Big( (e \pm f) \rho_1' - (e \mp f) \frac{\rho_1}{r} 
\Big) \, r \, dr = 0,
\end{equation}
we have
\begin{align*}
Q_{R, 1}(e) + Q_{R, - 1}(f) & + I_R(e \pm f)\\
& \geq \kappa_0 \Big( \| e \|_{H_1}^2 + \| f \|_{H_{- 1}}^2 + I_R(e\pm f) \Big) 
- C \int_R^{2 R} \frac{|e|^2 + |f|^2}{r^2} \, r \, dr.
\end{align*}
\end{prop}

\begin{proof}
We split the quantity $Q_{R, 1}(e) + Q_{R, - 1}(f) + I_R(e \pm f)$ as
$$
Q_\text{loc}^\pm \big( \chi_R e,\chi_R f \big) + Q_{\infty} \big( (1 - \chi_R) e, (1 
- \chi_R) f \big) + \boR_R \big( e, f \big),
$$
where
$$
Q_\text{loc}^\pm(u, v) = \int_0^{+ \infty} \Big( |u'|^2 + |v'|^2 + \frac{4}{r^2} 
|u|^2 - (1 - \rho_1^2) \big( |u|^2 + |v|^2 \big) + \rho_1^2 |u \pm v|^2 \Big) \, 
r \, dr,
$$
$$
Q_\infty(u, v) := \int_0^{+ \infty} \Big( |u'|^2 + |v'|^2 + \Big( \frac{4 - 
2\rho_1^2}{r^2} - (1 - \rho_1^2) \Big) |u|^2 + \Big( \frac{2 \rho_1^2}{r^2} - (1 
- \rho_1^2) \Big) |v|^2 \Big) \, r \, dr,
$$
and
\begin{align*}
\boR_R(u, v) := & \int_0^{+ \infty} 2 \Big( (u \chi_R)' (u (1 - \chi_R) )' + (v 
\chi_R)' (v (1 - \chi_R))' \Big) \, r \, dr \\ & + \int_0^{+ \infty} 2 \chi_R (1 
- \chi_R) \Big( \frac{4}{r^2} |u|^2 - (1 - \rho_1^2) \big( |u|^2 + |v|^2 \big) 
\Big) \, r \, dr.
\end{align*}
We now control each of the previous quantity separately.

\begin{step}
\label{S1}
If $u$ and $v$ are supported outside of the interval $[0, R]$, and $R_0$ is 
larger than some universal constant, then we deduce from the decay properties of 
$1 - \rho_1^2$ in Lemma~\ref{lem:propprofil} that
\begin{equation}
\label{eq:estimQinf}
Q_\infty(u, v) \geq \int_0^{+ \infty} \Big( |u'|^2 + |v'|^2 + \frac{1}{2r^2} 
\big( |u|^2 + |v|^2 \big) \Big) \, r \, dr \geq \kappa_0 \big( \| u \|_{H_1}^2 + 
\| v \|_{H_{- 1}}^2 \big).
\end{equation}
\end{step}

\begin{step}
\label{S2} 
We next claim that there exists a universal constant $C > 0$ such that
\begin{equation}
\label{eq:reste}
\boR_R(u, v) \geq - C \int_R^{2 R} \frac{|u|^2 + |v|^2}{r^2} \, r \, dr.
\end{equation}
for any $u \in H_1$ and any $v \in H_{- 1}$. Concerning the second integral in 
the definition of $\boR_R(u, v)$ , this estimate follows from the fact that 
$\chi_R (1 - \chi_R)$ is supported in $[R, 2R]$ and from the decay properties of 
$1 - \rho_1^2$ in Lemma~\ref{lem:propprofil}. For the first one, an integration 
by parts provides
\begin{align*}
\int_0^{+ \infty} 2 (u \chi_R)' \big( u (1 & - \chi_R) \big)' \, r \, dr = 2 
\int_R^{2 R} \chi_R (1 - \chi_R) |u'|^2 \, r \, dr \\
& - \int_R^{2 R} \Big( 2(\chi_R')^2 + (\chi_R - \chi_R^2)'' + \frac{(\chi_R - 
\chi_R^2)'}{r} \Big) |u|^2 \, r \, dr,
\end{align*}
where the first term is non-negative and
$$
\Big| 2 (\chi_R')^2 + (\chi_R - \chi_R^2)'' + \frac{(\chi_R - \chi_R^2)'}{r} 
\Big| \leq \frac{C}{R^2},
$$
pointwise on $[R, 2R]$. The analogous inequality holds for $v$, and 
inequality~\eqref{eq:reste} therefore follows.
\end{step}

\begin{step}
\label{S3}
We finally claim that, if $u \in H_1$ and $v \in H_{-1}$ are real-valued, 
compactly supported, and satisfy the orthogonality condition
\begin{equation}
\label{eq:orthoici}
\int_0^{+ \infty} \chi \Big( (u \pm v) \rho_1' - (u \mp v) \frac{\rho_1}{r} 
\Big) \, r \, dr = 0,
\end{equation}
then
\begin{equation}
\label{eq:estimQloc}
Q_\text{loc}^\pm(u, v) \geq \kappa_0 \bigg( \| u \|_{H_1}^2 + \| v \|_{H_{- 1}}^2 + 
\int_0^{+ \infty} \rho_1^2 |u \pm v|^2 \, r \, dr \bigg).
\end{equation}
In order to prove this claim, we first recall that it was proved 
in~\cite{dePiFeK1} that the quantities $Q_{\rm loc}^\pm(u, v)$ are non-negative 
and vanish if and only if $u \pm v = c \rho_1'$ and $u \mp v = - c \rho_1/r$ for 
some constant $c \in \R$. As a consequence, the quadratic forms
$$
\boQ_\text{loc}^\pm(u, v) := Q_\text{loc}^\pm(u, v) + \bigg( \int_0^{+ \infty} \chi \Big( 
(u \pm v) \rho_1' - (u \mp v) \frac{\rho_1}{r} \Big) \, r \, dr \bigg)^2,
$$
are positive definite on the Hilbert spaces $G^\pm$ associated to the norm 
$$
\big\| (u, v) \big \|_{G^\pm}^2 := \| u \|_{H_1}^2 + \| v \|_{H_{- 1}}^2 + \int_0^{+ 
\infty} \rho_1^2 |u \pm v|^2 \, r \, dr.
$$
The conclusion will follow from Lemma~\ref{lem:fred} if we can check that 
the quadratic forms $\boQ_\text{loc}^\pm$ are associated to a compact 
perturbation of the identity in $G^\pm$ for some equivalent norm.

We choose the norm given by 
$$
\big| \big| \big| (u, v) \big| \big| \big|_{G^\pm}^2 := \int_0^{+ \infty} \Big( 
|u'|^2 + |v'|^2 + \frac{4}{r^2} |u|^2 + |u \pm v|^2 \Big) \, r \, dr - 
\int_4^\infty \frac{|u|^2 + |v|^2}{r^2} \, r \, dr.
$$
We can readily check that $||| (u, v) |||_{G^\pm} \leq C \| (u, v) \|_{G^\pm}$ 
for some universal constant $C >0$. The converse inequality requires some 
explanation, actually even the fact that $||| \cdot |||_{G^\pm}$ defines a norm. 
For that purpose, computing the discriminant gives the pointwise inequality
$$
\Big( \frac{2}{r^2} + \frac{1}{2} \Big) |u|^2 \pm u v + \Big( \frac{1}{2} - 
\frac{3}{2 r^2} \Big) |v|^2 > 0,
$$
for $(u, v) \neq (0, 0)$ and $r > \sqrt{12}$. Therefore, we obtain 
$$
\frac{4}{r^2} |u|^2 + |u \pm v|^2 - \frac{1}{r^2} \big( |u|^2 + |v|^2 \big) \geq 
\frac{1}{r^2} |u|^2 + \frac{1}{2} \rho_1^2 |u \pm v|^2 + \frac{1}{2 r^2} |v|^2,
$$
for $r \geq 4$. On the other hand, there exists some universal constant $C > 0$ 
such that
$$
\frac{4}{r^2} |u|^2 + |u \pm v|^2 \geq \frac{1}{C} \Big( \frac{4}{r^2} |u|^2 + 
|v|^2 + \rho_1^2 |u \pm v|^2 \Big),
$$
when $0 \leq r \leq 4$. Combining the previous estimates is enough to guarantee 
that $||| \cdot |||_{G^\pm}$ is a norm and is equivalent to $\| \cdot 
\|_{G^\pm}$. Next, we write
\begin{align*}
Q_\text{loc}^\pm(u, v) - \big| \big| \big| (u, v) \big| \big| \big|_{G^\pm}^2 &
= \int_4^{+ \infty} \Big( \frac{1}{r^2} - (1 -\rho_1^2) \Big) \big( |u|^2 + 
|v|^2 \big) \, r \, dr \\
& - \int_0^4 (1 - \rho_1^2) \big( u^2 + v^2 \big) \, r \, dr - \int_0^{+ \infty} 
(1 - \rho_1^2) \big| u \pm v \big|^2 \, r \, dr,
\end{align*}
and each of the three terms on the right-hand-side is compact in view of the 
decay properties of $\rho_1$. Finally, $\boQ_\text{loc}^\pm - Q_\text{loc}^\pm$ 
is also compact as the square of a scalar product.
\end{step} 

We are now in position to complete the proof of Proposition~\ref{prop:H1Hm1}. We 
apply Step~\ref{S3} to $u = \chi_R e$ and $v = \chi_R f$. Since $\chi_R \equiv 
1$ on the interval $[0, 2]$, the orthogonality condition~\eqref{eq:ortho1m1} 
implies the orthogonality condition~\eqref{eq:orthoici}. Therefore, we have
\begin{equation}
\label{eq:step4_1}
Q_\text{loc}^\pm(\chi_R e, \chi_R f) \geq \kappa_0 \bigg( \| \chi_R e \|_{H_1}^2 
+ \| \chi_R f \|_{H_{- 1}}^2 + \int_0^{+ \infty} \chi_R^2 \rho_1^2 |e \pm f|^2 
\, r \, dr \bigg).
\end{equation}
Next, we apply Step~\ref{S1} with $u = (1 - \chi_R) e$ and $v = (1 - \chi_R) f$ 
so as to obtain
\begin{equation}
\label{eq:step4_2}
Q_\infty \big( (1 - \chi_R) e, (1 - \chi_R) f \big) \geq \kappa_0 \Big( \big\| 
(1 - \chi_R) e \big\|_{H_1}^2 + \big\| (1 - \chi_R) f \big\|_{H_{- 1}}^2 \Big).
\end{equation}
The summation of~\eqref{eq:step4_1} and~\eqref{eq:step4_2} combined with 
Step~\ref{S2} yields the conclusion.
\end{proof}

We finally conclude the

\begin{proof}[Proof of Proposition~\ref{prop:quad}]
In view of~\eqref{eq:nabla-V1}, we first note that the orthogonality conditions 
in~\eqref{eq:orthocond} translate into
$$
\int_{\R^2} \chi \langle \varepsilon, i V_1\rangle_\C = \int_0^{+ \infty} \chi 
b_0 \rho_1 \, r \, dr = 0,
$$
$$
\int_{\R^2} \chi \langle \varepsilon, \partial_x V_1 \rangle_\C = \frac{1}{2} 
\int_0^{+ \infty} \chi \Big( (a_1 + a_{- 1}) \rho_1' - (a_1 - a_{- 1}) 
\frac{\rho_1}{r} \Big) \, r \, dr = 0,
$$
and
$$
\int_{\R^2} \chi \langle \varepsilon, \partial_y V_1 \rangle_\C = - \frac{1}{2} 
\int_0^{+ \infty} \Big( (b_1 - b_{- 1}) \rho_1' - (b_1 + b_{- 1}) 
\frac{\rho_1}{r} \Big) \, r \, dr = 0.
$$
As a consequence, we can estimate the terms in~\eqref{eq:decompboQ} 
and~\eqref{eq:decompboI} using Corollaries~\ref{cor:niveau0},~\ref{cor:jlarger3} 
and~\ref{cor:niveau2}, and Proposition~\ref{prop:H1Hm1} in order to obtain the 
lower bound
$$
Q_R(\varepsilon) \geq \kappa_0 \bigg( \sum_{j \in \Z} \| \varepsilon_j 
\|_{H_j}^2 + \boI_R(\varepsilon) \bigg) - C \int_R^{2 R} \frac{|\varepsilon_1|^2 
+ |\varepsilon_{- 1}|^2}{r^2} \, r \, dr,
$$
for some constant $\kappa_0 > 0$ independent of $R$. Fix $N_0 \geq 1$. By the 
pigeon-hole principle, we can find $R_0 \leq R \leq 2^{N_0} R_0$ such that
$$
\int_R^{2 R} \frac{|\varepsilon_1|^2 + |\varepsilon_{-1}|^2}{r^2} \, r \, dr 
\leq \frac{1}{N_0} \int_{R_0}^{2^{N_0}R_0} \frac{|\varepsilon_1|^2 + 
|\varepsilon_{-1}|^2}{r^2} \, r \, dr \leq \frac{2}{N_0} \Big( \| \varepsilon_1 
\|_{H_1}^2 + \| \varepsilon_{- 1} \|_{H_{- 1}}^2 \Big).
$$
We therefore choose $N_0$ such that $2 C/N_0 \leq \kappa_0/2$ and the conclusion 
follows with $\kappa_0$ replaced by $\kappa_0/2$.
\end{proof}

\section{Modulation parameters - Proof of Proposition~\ref{prop:modul-param}} 
\label{sec:proof-modul}

The proof of Proposition~\ref{prop:modul-param} is classical 
(see e.g.~\cite{BetGrSm1} and the references therein). 
However, we have to handle with care the norm $\| \cdot \|_H$ since it 
is not left invariant by translation. This is the reason why we 
provide the following detail.

The main ingredient is to apply the implicit function theorem to the map
\begin{equation}
\label{eq:def-Xi}
\Xi(\Psi, b, \varphi) = \bigg( \int_{\R^2} \chi \langle \varepsilon, \partial_x 
V_1 \rangle_\C, \int_{\R^2} \chi \langle \varepsilon, \partial_y V_1 \rangle_\C, 
\int_{\R^2} \chi \langle \varepsilon, i V_1 \rangle_\C \bigg).
\end{equation}
The function $\Psi$ in this expression belongs to $H$, the vector $b$ is in 
$\R^2$ and the number $\varphi$ in $\R$. As above, we have set $\varepsilon = 
e^{- i \varphi} \Psi(\cdot + b) - V_1$. The map $\Xi$ is well-defined from $H 
\times \R^2 \times \R$ to $\R^3$ and it satisfies
\begin{equation}
\label{eq:Xi-zero}
\Xi(e^{i \varphi} V_1(\cdot - b), b, \varphi) = 0.
\end{equation}
Applying the implicit function theorem, we can expect to construct parameters 
$a(\Psi) \in \R^2$ and $\varphi(\Psi) \in \R$ so as to guarantee the 
orthogonality conditions in~\eqref{eq:orthocond} for any function $\Psi$ in a 
neighbourhood of any vortex solution $e^{i \varphi} V_1(\cdot - b)$. We first 
perform this construction for the original vortex solution $V_1$.

\begin{lem}
\label{lem:implicit-V1}
Set $\bB_{V_1}(r) := \{ \Psi \in H \text{ s.t. } \| \Psi - V_1 \|_H < r \big\}$ 
for any $r > 0$. There exist $\rho > 0$ and $\Lambda > 0$ such that the 
following statements hold. There exist two maps $b := (b_1, b_2) \in 
\boC^1(\bB_{V_1}(\rho), \R^2)$ and $\varphi \in \boC^1(\bB_{V_1}(\rho), \R)$ 
such that, given any function $\Psi \in \bB_{V_1}(\rho)$, the pair $(b, \varphi) 
= (b(\Psi), \varphi(\Psi))$ is the unique solution in the product set $(- 
\Lambda \rho, \Lambda \rho)^3$ of the equations
$$
\int_{\R^2} \chi \langle \varepsilon, \partial_x V_1 \rangle_\C = \int_{\R^2} 
\chi \langle \varepsilon, \partial_y V_1 \rangle_\C = \int_{\R^2} \chi \langle 
\varepsilon, i V_1 \rangle_\C = 0,
$$
where $\varepsilon = e^{- i \varphi} \Psi(\cdot + b) - V_1$. Moreover, the maps 
$b$ and $\varphi$ satisfy
\begin{equation}
\label{eq:cont-Lambda}
|b_1(\Psi_2) - b_1(\Psi_1)| + |b_2(\Psi_2) - b_2(\Psi_1)| + |\varphi(\Psi_2) - 
\varphi(\Psi_1)| \leq \Lambda \| \Psi_2 - \Psi_1 \|_H,
\end{equation}
for any functions $(\Psi_1, \Psi_2) \in \bB_{V_1}(\rho)^2$.
\end{lem}

\begin{proof}
Since the map $\Xi$ is continuously differentiable from $H \times \R^2 \times 
\R$ to $\R^3$, we can compute
\begin{equation}
\label{eq:dXi-b}
\begin{split}
\nabla_b \Xi(\Psi, b, \varphi) = \bigg( & \int_{\R^2} \chi \langle e^{- i 
\varphi} \nabla \Psi(\cdot + b), \partial_x V_1 \rangle_\C, \int_{\R^2} \chi 
\langle e^{- i \varphi} \nabla \Psi(\cdot + b), \partial_y V_1 \rangle_\C,\\
& \int_{\R^2} \chi \langle e^{- i \varphi} \nabla \Psi(\cdot + b), i V_1 
\rangle_\C \Big),
\end{split}
\end{equation}
and
\begin{equation}
\label{eq:dXi-phi}
\begin{split}
\partial_\varphi \Xi(\Psi, b, \varphi) = - \bigg( & \int_{\R^2} \chi \langle i 
e^{- i \varphi} \Psi(\cdot + b), \partial_x V_1 \rangle_\C, \int_{\R^2} \chi 
\langle i e^{- i \varphi} \Psi(\cdot + b), \partial_y V_1 \rangle_\C,\\
& \int_{\R^2} \chi \langle i e^{- i \varphi} \Psi(\cdot + b), i V_1 \rangle_\C 
\bigg).
\end{split}
\end{equation}
Since
$$
\int_{\R^2} \chi \langle \partial_x V_1, \partial_y V_1 \rangle_\C = \int_{\R^2} 
\chi \langle i V_1, \partial_x V_1 \rangle_\C = \int_{\R^2} \chi \langle i V_1, 
\partial_y V_1 \rangle_\C = 0,
$$
we deduce from the previous formulae that the differential
\begin{equation}
\label{eq:def-D}
d_{b_1, b_2, \varphi} \Xi(V_1, 0, 0) = \begin{pmatrix} \int_{\R^2} \chi 
|\partial_x V_1|^2 & 0 & 0\\ 0 & \int_{\R^2} \chi |\partial_y V_1|^2 & 0\\ 0 & 0 
& \int_{\R^2} \chi |V_1|^2 \end{pmatrix},
\end{equation}
is a continuous isomorphism from $\R^2 \times \R$ to $\R^3$. In view 
of~\eqref{eq:Xi-zero}, we infer from the implicit function theorem the existence 
of some $\rho > 0$, of an open neighbourhood $\boU$ of $(V_1, 0 ,0)$ in 
$H \times \R^3$ and of two functions $b \in \boC^1(\bB_{V_1}(\rho), \R^2)$ and 
$\varphi \in \boC^1(\bB_{V_1}(\rho), \R)$ such that, for any datum $(\Psi, b, 
\varphi) \in \boU$, the equation $\Xi(\Psi, b, \varphi) =0$ owns a unique 
solution given by $(b, \varphi) = (b(\Psi), \varphi(\Psi))$. By continuous 
differentiability of the map $\gamma := (b, \varphi)$, we can decrease the value 
of $\rho$ such that the operator norm $\| d\gamma(\Psi) \|$ of the differentials 
$d\gamma(\Psi)$ is less than $\Lambda := 1 + \| d\gamma(V_1) \|$ on the ball 
$\bB_{V_1}(\rho)$. Inequality~\eqref{eq:cont-Lambda} then follows from the mean 
value inequality. In turn, we infer from~\eqref{eq:cont-Lambda} that the map 
$\gamma$ is valued into the ball $B_{\Lambda \rho}$ and we can decrease the 
value of $\rho$, if necessary, so as to replace the open subset $\boU$ by the 
product set $\bB_{V_1}(\rho) \times B_{\Lambda \rho}$. This ends the proof of 
Lemma~\ref{lem:implicit-V1}.
\end{proof}

We now extend the previous construction to the neighbourhood of any fixed vortex 
solution $e^{i \phi} V_1(\cdot - a)$ by using the translation and phase 
invariances.

\begin{cor}
\label{cor:implicit-V1}
For $a \in \R^2$ and $\phi \in \R$, consider the balls $\bB_{(a, \phi)}(r) := 
\big\{ \Psi \in H\text{ s.t. }\| e^{- i \phi} \Psi(\cdot \linebreak[0] + a) - 
V_1 \|_H < r \big\}$ for any $r > 0$, and set
$$
b_{a, \phi}(\psi) = a + b \big( e^{- i \phi} \psi(\cdot + a) \big) \quad 
\text{and} \quad \varphi_{a, \phi}(\psi) = a + \varphi \big( e^{- i \phi} 
\psi(\cdot + a) \big),
$$
for any function $\psi \in \bB_{(a, \phi)}(\rho)$. Given any function $\Psi \in 
\bB_{(a, \phi)}(\rho)$, the pair $(b, \varphi) = (b_{a, \phi}(\Psi), 
\linebreak[0] \varphi_{a, \phi}(\Psi))$ is the unique solution in the product 
set $(a_1 -\Lambda \rho, a_1 + \Lambda \rho) \times (a_2 -\Lambda \rho, a_2 + 
\Lambda \rho) \times (\phi -\Lambda \rho, \phi + \Lambda \rho)$ of the equations
$$
\int_{\R^2} \chi \langle \varepsilon, \partial_x V_1 \rangle_\C = \int_{\R^2} 
\chi \langle \varepsilon, \partial_y V_1 \rangle_\C = \int_{\R^2} \chi \langle 
\varepsilon, i V_1 \rangle_\C = 0,
$$
where $\varepsilon = e^{- i \varphi} \Psi(\cdot + b) - V_1$. Moreover, the maps 
$b_{a, \phi}$ and $\varphi_{a, \phi}$ satisfy
\begin{align*}
\big| [b_{a, \phi}]_1(\Psi_2) - [b_{a, \phi}]_1(\Psi_1) \big| + \big| [b_{a, 
\phi}]_2(\Psi_2) - [b_{a, \phi}]_2(\Psi_1) \big| & + \big| \varphi_{a, 
\phi}(\Psi_2) - \varphi_{a, \phi}(\Psi_1) \big|\\
\leq & \Lambda \big\| \Psi_2(\cdot + a) - \Psi_1(\cdot + a) \big\|_H,
\end{align*}
for any functions $(\Psi_1, \Psi_2) \in \bB_{(a, \phi)}(\rho)^2$.
\end{cor}

\begin{proof}
Corollary~\ref{cor:implicit-V1} is a direct consequence of 
Lemma~\ref{lem:implicit-V1} once we have observed that the map $(\Psi, b, \phi) 
\mapsto (e^{- i \phi} \Psi(\cdot + a), b - a, \varphi - \phi)$ is a bijection 
from the product set $\bB_{(a, \phi)}(\rho) \times (a_1 - \Lambda \rho, a_1 + \Lambda \rho) \times 
(a_2 - \Lambda \rho, a_2 + \Lambda \rho) \times (\phi - \Lambda \rho, \phi + 
\Lambda \rho)$ onto $\bB_{V_1}(\rho) \times (- \Lambda \rho, \Lambda \rho)^3$.
\end{proof}

The next step in the proof of Proposition~\ref{prop:modul-param} is to extend 
the previous construction to a neighbourhood of the orbit of $V_1$ of the form 
$\boV(\alpha)$. By definition, this neighbourhood is equal to
$$
\boV(\alpha) = \underset{(a, \phi) \in \R^3}{\cup} \bB_{(a, \phi)}(\alpha).
$$
For $\alpha \leq \rho$, the existence of modulation parameters so that the 
orthogonality conditions in~\eqref{eq:orthocond} are satisfied results from 
Corollary~\ref{cor:implicit-V1}. In order to complete the proof of 
Proposition~\ref{prop:modul-param}, it essentially remains to establish that the 
choice of these parameters can be made in a continuously differentiable way.

In this direction, the main difficulty is to prevent the possibility that a 
function $\Psi \in \boV(\rho)$ belongs to two balls $\bB_{(a_1, \phi_1)}(\rho)$ and 
$\bB_{(a_2, \phi_2)}(\rho)$ for points $a_1$ and $a_2$ at a large distance from one 
another. In this case, the translation parameter $a$ can be chosen either close 
to $a_1$ or to $a_2$ in view of~\eqref{eq:cont-Lambda}. Hence it is not so 
direct to find a continuously differentiable choice for this parameter.

In order to by-pass this difficulty, we first show the following lemma.

\begin{lem}
\label{lem:ball-close}
There exist $\alpha_0 > 0$ and $R_0 >0$ such that, if 
$$\bB_{(a_1, \phi_1)}(\alpha_0) \cap \bB_{(a_2, \phi_2)}(\alpha_0) \neq \emptyset,$$
for points $(a_1, a_2) \in \R^4$ and numbers $(\phi_1, \phi_2) \in\R^2$, then
$$|a_1 - a_2| \leq R_0.$$
\end{lem}

\begin{proof}
Let $\alpha > 0$ and consider a function $\Psi \in \bB_{(a_1, \phi_1)}(\alpha) \cap 
\bB_{(a_2, \phi_2)}(\alpha)$ for two points $a_1$ and $a_2$, and two numbers 
$\phi_1$ and $\phi_2$. By definition, this function first satisfies
\begin{equation}
\label{haouas}
\int_{\R^2} \big( 1 - |V_1|^2 \big) \Big| \nabla \big( \Psi(\cdot + a_1) - e^{i 
\phi_1} V_1 \big) \Big|^2 < \alpha^2.
\end{equation}
In view of Lemma~\ref{lem:propprofil}, the integral
$$
I_1 := \big\| V_1 \big\|_H^2 = \int_\R \Big( \big| \nabla |V_1|^2 \big|^2 + (1 - 
|V_1|^2) \big| \nabla V_1 \big|^2 \Big),
$$
is finite and positive. Moreover, we can find $R > 0$ such that
\begin{equation}
\label{marchand}
\int_{B_R} \Big( \big| \nabla |V_1|^2 \big|^2 + (1 - |V_1|^2) \big| \nabla V_1 
\big|^2 \Big) = \frac{31 I_1}{32}.
\end{equation}
Hence we infer from~\eqref{haouas} and the inequality $(\alpha - \beta)^2 
\geq \alpha^2/2 - \beta^2$ that
$$
\int_{B_R} (1 - |V_1|^2) \big| \nabla \Psi(\cdot + a_1) \big|^2 \geq \frac{31 
I_1}{32} - \alpha^2.
$$
Since $\rho_1 = |V_1| < 1$, we are led to
\begin{equation}
\label{baille}
\int_{B_R(a_1)} \big| \nabla \Psi \big|^2 \geq \frac{31 I_1}{64} - \alpha^2.
\end{equation}
Here we have set, as in the sequel, $B_r(a) = \{ x \in \R^2 \text{ s.t.} |x - a| 
< r \}$ for any $r > 0$ and any $a \in \R^2$. 

Similarly, we know that
\begin{equation}
\label{taofifenua}
\int_{\R^2} \Big| \nabla \big( \bar{V}_1 (\Psi(\cdot + a_2) - e^{i \phi_2} V_1) 
\big) \Big|^2 < \alpha^2,
\end{equation}
so that, by~\eqref{marchand},
$$
\int_{B_R(a_2)^c} \Big| \nabla \big( \bar{V}_1(\cdot - a_2) \Psi \big) \Big|^2 
\leq 2 \Big( \frac{I_1}{32} + \alpha^2 \Big).
$$
Assuming that $|a_2 - a_1| \geq 2 R$, we obtain
$$
\frac{1}{2} \int_{B_R(a_1)} |V_1(\cdot - a_2)|^2 \big| \nabla \Psi \big|^2 \leq 
2 \Big( \frac{I_1}{32} + \alpha^2 \Big) + \int_{B_R(a_1)} |\Psi|^2 \big| \nabla 
V_1(\cdot - a_2) \big|^2.
$$
At this stage, we can increase, if necessary, the value of the number $R$ 
such that $|V_1(x)|^2 \geq 1/2$ for $|x| \geq R$. In this case, we obtain
$$
\int_{B_R(a_1)} \big| \nabla \Psi \big|^2 \leq 8 \Big( \frac{I_1}{32} + \alpha^2 
\Big) + 4 \int_{B_R(a_1)} |\Psi|^2 \big| \nabla V_1(\cdot - a_2) \big|^2.
$$
In view of the existence of a universal constant $C > 0$ such that $|\nabla 
V_1(x)|^2 \leq C/(1 + |x|^2)$, we next have
$$
\int_{B_R(a_1)} \big| \nabla \Psi \big|^2 \leq 8 \Big( \frac{I_1}{32} + \alpha^2 
\Big) + \frac{4 C}{1 + (|a_2 - a_1| - R)^2} \int_{B_R(a_1)} |\Psi|^2.
$$
Going to the proof of Lemma~\ref{lem:comp-emb}, we next find $K_R > 0$, 
depending only on $R$, such that
$$
\int_{B_R(a_1)} |\Psi|^2 \leq K_R \big\| \Psi \big\|_{H_{a_1}}^2 \leq 2 K_R 
\Big( \big\| \Psi - e^{i \varphi_1} V_1(\cdot - a_1) \big\|_{H_{a_1}}^2 + I_1 
\big),
$$
for any number $\varphi \in \R$. Hence we are led to
\begin{equation}
\label{willemse}
\int_{B_R(a_1)} \big| \nabla \Psi \big|^2 \leq 8 \Big( \frac{I_1}{32} + \alpha^2 
\Big) + \frac{8 C K_R}{1 + (|a_2 - a_1| - R)^2} \Big( \alpha^2 + I_1 \Big).
\end{equation}

Combining~\eqref{baille} and~\eqref{willemse} next gives
$$
\frac{31 I_1}{64} - \alpha^2 \leq 8 \Big( \frac{I_1}{32} + \alpha^2 \Big) 
+ \frac{8 C K_R}{1 + (|a_2 - a_1| - R)^2} \Big( \alpha^2 + I_1 \Big),
$$
and we can choose $\alpha = \alpha_0 := \sqrt{I_1/64}$ in order to obtain
$$
|a_2 - a_1| \leq R + 10 \sqrt{C K_R}.
$$
The conclusion follows for $R_0 = \max \{ 2 R, R + 10 \sqrt{C K_R} \}$.
\end{proof}

We next refine the bound in Lemma~\ref{lem:ball-close} for $\alpha$ small.

\begin{lem}
\label{lem:ball-very-close}
Let $\mu > 0$. There exists $\nu > 0$ such that, if
$$\bB_{(a_1, \phi_1)}(\nu) \cap \bB_{(a_2, \phi_2)}(\nu) \neq \emptyset,$$
for points $(a_1, a_2) \in \R^4$ and numbers $(\phi_1, \phi_2) \in\R^2$, then
$$\big| a_2 - a_1 \big| + \big| e^{i \phi_2} - e^{i \phi_1} \big| < \mu.$$
\end{lem}

\begin{proof}
We argue by contradiction assuming that the statement in 
Lemma~\ref{lem:ball-very-close} is wrong. In this case, we can find $\mu > 0$, 
as well as sequences $(a_n)_{n \in \N}$, $(b_n)_{n \in \N}$, $(\phi_n)_{n \in 
\N}$ and $(\varphi_n)_{n \in \N}$ such that
\begin{equation}
\label{leroux}
\bB_{(a_n, \phi_n)} \Big( \frac{1}{2^n} \Big) \cap \bB_{(b_n, \varphi_n)} \Big( 
\frac{1}{2^n} \Big) \neq \emptyset,
\end{equation}
and
\begin{equation}
\label{rebbadj}
\big| b_n - a_n \big| + \big| e^{i \varphi_n} - e^{i \phi_n} \big| \geq \mu,
\end{equation}
for any $n \in \N$. Up to a subsequence, there exist two numbers $\phi_1$ 
and $\phi_2$ such that
\begin{equation}
\label{gros}
e^{i \phi_n} \to e^{i \phi_\infty} \quad \text{and} \quad e^{i \varphi_n} 
\to e^{i \varphi_\infty},
\end{equation}
as $n \to + \infty$. For $n$ large enough, we also infer from~\eqref{leroux} 
that the difference $a_n - b_n$ is bounded by the number $R_0$ in 
Lemma~\ref{lem:ball-close}. Up to a further subsequence, we can assume the 
existence of a point $d_\infty \in \R^2$ such that
\begin{equation}
\label{pesenti}
b_n - a_n \to d_\infty,
\end{equation}
as $n \to + \infty$. Since the norm $|d_\infty|$ is positive by~\eqref{rebbadj}, 
it follows that
\begin{equation}
\label{cazeaux}
\big| b_n - a_n \big| \leq 2 |d_\infty|,
\end{equation}
for $n$ large enough. 

Consider next a sequence of functions $\Psi_n$ in $\bB_{(a_n, \phi_n)}(1/2^n) 
\cap \bB_{(b_n, \varphi_n)}(1/2^n)$. Going back to the proof of 
Lemma~\ref{lem:comp-emb}, we can find $C > 0$, depending only on $|d_\infty|$, 
such that
$$
\int_{B_{3 |d_\infty|}} \Big( \big| e^{- i \phi_n} \Psi_n(\cdot + a_n) - V_1 
\big|^2 + \big| e^{- i \varphi_n} \Psi_n(\cdot + b_n) - V_1 \big|^2 \Big) \leq 
\frac{C}{2^n},
$$
for any $n \in \N$. In view of~\eqref{cazeaux}, we observe that
\begin{align*}
\int_{B_{|d_\infty|}} \big| e^{i (\varphi_n - \phi_n)} V_1(\cdot & - b_n + a_n) 
- V_1 \big|^2 \leq 2 \int_{B_{3 |d_\infty|}} \big| e^{- i \phi_n} \Psi_n(\cdot + 
a_n) - V_1 \big|^2 \\
& + 2 \int_{B_{3 |d_\infty|}(b_n - a_n)} \big| e^{- i \varphi_n} \Psi_n(\cdot + 
a_n) - V_1(\cdot - b_n + a_n) \big|^2,
\end{align*}
so that
$$
\int_{B_{|d_\infty|}} \big| e^{i (\varphi_n - \phi_n)} V_1(\cdot - b_n + a_n) - V_1 \big|^2 \leq \frac{4 C}{2^n}.
$$
Combining with~\eqref{gros}and~\eqref{pesenti}, we are led to
$$
\int_{B_{|d_\infty|}} \big| e^{i (\phi_2 - \phi_1)} V_1(\cdot - d_\infty) - V_1 \big|^2 = 0,
$$
in the limit $n \to + \infty$. We conclude that $e^{i (\phi_2 - \phi_1)} V_1(\cdot - d_\infty) = V_1$ 
on the ball $B_{|d_\infty|}$. Since the function $V_1$ only vanishes at the origin, 
we infer that $d_\infty = 0$, and then $e^{i \phi_2} = e^{i \phi_1}$. 
This is a contradiction with the inequality $|d_\infty| + |e^{i \phi_2} - e^{i \phi_1}| \geq \mu$, 
that follows from~\eqref{rebbadj} in the limit $n \to + \infty$. This completes the proof 
of Lemma~\ref{lem:ball-very-close}.
\end{proof}

We are now in position to conclude the

\begin{proof}[Proof of Proposition~\ref{prop:modul-param}]
 Without loss of generality, we can assume that the numbers $\rho$ and $\Lambda$ 
 in Lemma~\ref{lem:implicit-V1} and Corollary~\ref{cor:implicit-V1} satisfy $\rho \Lambda < 1$. 
 We consider the number $\nu$ provided by Lemma~\ref{lem:ball-very-close} for 
 $\mu = \rho \Lambda/8$ and we set $\alpha := \min \{ \rho/2, \nu \}$.

Given a function $\Psi \in \boV(\alpha)$, we can find a point $a$ and a number
$\phi$ such that $\Psi$ belongs to the ball $\bB_{(a, \phi)}(\alpha)$. Since $\alpha
< \rho$, we infer from Corollary~\ref{cor:implicit-V1} the existence of a point
$a(\Psi) = b_{a, \phi}(\Psi) \in \R^2$ and of a number $\varphi(\Psi) =
\varphi_{a, \phi}(\Psi) \in \R$ such that the orthogonality conditions
in~\eqref{eq:orthocond} are satisfied. 
We claim that the value of $a(\Psi)$ does not depend on the choice of $a$ and $\phi$. 
The number $\varphi(\Psi)$ is also independent of this choice, but modulo $2 \pi$.

Assume indeed that the function $\Psi$ is in $\bB_{(b, \vartheta)}(\alpha)$ for another 
point $b \in \R^2$ and another number $\vartheta \in \R$. 
The intersection $\bB_{(a, \phi)}(\nu) \cap \bB_{(b, \vartheta)}(\nu)$ is then not empty, 
so that by Lemma~\ref{lem:ball-very-close},
\begin{equation}
\label{ollivon}
\big| b - a \big| + \big| e^{i \vartheta} - e^{i \phi} \big| < \mu = \frac{\rho \Lambda}{8}.
\end{equation}
Recall here that $|e^{i t} - 1| = 2 |\sin(t/2)| \geq 2 |t|/\pi$ when $t 
\in [-
\pi, \pi]$. Hence, there exists an integer $k \in \Z$ such that
\begin{equation}
\label{cretin}
|\vartheta + 2 \pi k - \phi| < \frac{\pi \mu}{2} < \frac{\rho \Lambda}{4}.
\end{equation}
On the other hand, we also infer from Corollary~\ref{cor:implicit-V1} that
$$
\big| a_1(\Psi) - a_1 \big| + \big| a_2(\Psi) - a_2 \big| + \big| \varphi(\Psi) 
- \phi \big| \leq \Lambda \big\| \Psi(\cdot + a) - e^{i \varphi} V_1(\cdot + a - a) \|_H < \frac{\Lambda \rho}{2}.
$$
Combining with~\eqref{ollivon} and~\eqref{cretin}, we obtain
$$
\big| a_1(\Psi) - b_1 \big| + \big| a_2(\Psi) - b_2 \big| + \big| \varphi(\Psi) 
- \vartheta - 2 \pi k \big| < \Lambda \rho.
$$
Since $\Xi(\Psi, a(\Psi), \varphi(\Psi) + 2 k \pi) = 0$, we deduce from 
Corollary~\ref{cor:implicit-V1} that $a(\Psi) = b_{b, \vartheta}(\Psi)$ 
and $\varphi(\Psi) = \varphi_{b, \vartheta}(\Psi) + 2 k \pi$. In conclusion, 
the choice of $a(\Psi)$ and $\varphi(\Psi)$ (modulo $2 \pi$) does not depend 
on the choice of $a$ and $\phi$ such that $\Psi \in \bB_{(a, \phi)}(\alpha)$. 
Therefore, the maps $a$ and $\varphi$ are well-defined from $\boV(\alpha)$ 
with values in $\R^2 \times \R / 2 \pi \Z$. They are continuous differentiable 
on $\boV(\alpha)$ due to the continuous differentiability of the maps 
$b_{a, \phi}$ and $\varphi_{a, \phi}$.

We now turn to the proof of~\eqref{eq:est-modul}. 
When $\big\| e^{- i \vartheta} \Psi(\cdot + b) - V_1 \big\|_H < \alpha$, 
we have $a(\Psi) = b_{b, \vartheta}(\Psi)$ and $\varphi(\Psi) = \varphi_{b, \vartheta}(\Psi)$. 
As a consequence of Corollary~\ref{cor:implicit-V1}, we first obtain
\begin{equation}
\label{alldritt}
\big| a(\Psi) - b \big| + \big| \varphi(\Psi) - \vartheta \big| \leq 
\sqrt{2} \Lambda \big\| e^{- i \vartheta}\Psi(\cdot + b) - V_1 \|_H.
\end{equation}
On the other hand, we know that
$$
\big\| \varepsilon \big\|_H = \big\| e^{- i \varphi(\Psi)} \Psi(\cdot + 
a(\Psi)) - V_1 \big\|_H.$$
Since $|a(\Psi) - b| < \sqrt{2} \Lambda \alpha < 1/\sqrt{2}$ by~\eqref{alldritt}, 
we can invoke the uniform boundedness of the translation operators in 
Lemma~\ref{lem:bdd-translation} in order to find $C > 0$ such that
$$
\big\| \varepsilon \big\|_H \leq C \big\| e^{- i \varphi(\Psi)} \Psi(\cdot + b) - V_1(\cdot - a(\Psi) +b) \big\|_H.
$$
In particular, we obtain
$$
\big\| \varepsilon \big\|_H \leq C \Big( \big\| e^{- i \vartheta} \Psi(\cdot + b) - V_1 \big\|_H + \big| e^{i (\vartheta - \varphi(\Psi))} - 1 \big| \| V_1 \|_H + \big\| V_1 - V_1(\cdot - a(\Psi) +b) \big\|_H \Big).
$$
Again since $|a(\Psi) - b| < 1/\sqrt{2}$ by~\eqref{alldritt}, we can infer 
from Lemma~\ref{lem:Lipschitz-translation} the existence of a further $C > 0$ such that
$$
\big\| \varepsilon \big\|_H \leq C \Big( \big\| e^{- i \vartheta} \Psi(\cdot + b) - V_1 \big\|_H + \big| \varphi(\Psi) - \vartheta \big| \| V_1 \|_H + \big| a(\Psi) - b \big| \Big).
$$
Estimate~\eqref{eq:est-modul} then follows from~\eqref{alldritt}. This completes 
the proof of Proposition~\ref{prop:modul-param}.
\end{proof}

\section{Evolution of the modulation parameters - Pr\-oof of 
Proposition~\ref{prop:evol-param}}

Let $\tau > 0$ to be fixed later. Under the assumption $d_E(V_1, \Psi_0) \leq 
\tau$, we can go back to the proofs of Theorems~\ref{thm:estimcore} 
and~\ref{thm:stab} in order to check that the solution $\Psi_t$ lies in a set 
$\boV(\alpha_\tau)$ for any $t \in \R$. Here, the numbers $\alpha_\tau$ 
tend to $0$ when $\tau \to 0$. In particular, we can apply 
Proposition~\ref{prop:modul-param} for $\tau$ small enough. This provides 
modulation parameters $a(t) := a(\Psi_t) \in \R^2$ and $\varphi(t) := 
\varphi(\Psi_t) \in \R/2 \pi \Z$ that satisfy all the statements in 
Proposition~\ref{prop:modul-param}.

Recall here that the solution $\Psi$ lies in $\boC(\R, \Psi_0 + H^1(\R^2)$ by Proposition~\ref{prop:cauchy}. Going back to the proof of Lemma~\ref{lem:Hilbert-H}, we check that it remains continuous with values in the Hilbert space $H$. As a consequence of Proposition~\ref{prop:modul-param}, the previous maps $a$ and $\varphi$ are also continuous from $\R$ to $\R^2$, respectively $\R/ 2 \pi \Z$. Up to the choice of a constant in $2 \pi \Z$, we can therefore reduce the map $\varphi$ to a continuous real-valued function.

In order to prove the continuous differentiability of $a$ and $\varphi$, we rely on the decomposition of any function in $E$ given by Lemma~\ref{lem:pourCauchy}. We first assume that the initial datum $\Psi_0$ takes the form $U_0 + w_0$, with $U_0 \in \boU$ and $w_0 \in H^3(\R^2)$. In this case, the corresponding solution $\Psi$ belongs to $\boC^0(\R, U_0 + H^3(\R^2))$ by Proposition~\ref{prop:cauchy2}. In view of~\eqref{eq:GP}, it is also in $\boC^1(\R, U_0 + H^3(\R^2))$, so in $\boC^1(\R, H)$ by Lemma~\ref{lem:Hilbert-H}. In view of Proposition~\ref{prop:modul-param}, the functions $a$ and $\varphi$ are then in $\boC^1(\R, \R^2)$, respectively $\boC^1(\R, \R)$. In order to extend this property to all initial data in $E$, we now compute the time derivatives $a'$ and $\varphi'$ by differentiating the orthogonality conditions in~\eqref{eq:orthocond}. We will eventually rely on these computations and a standard density argument in order to complete the proof of Proposition~\ref{prop:evol-param}.

Before going into to the computations of the derivatives $a'$ and $\varphi'$, we deduce from the previous smoothness properties of the maps $a$ and $\varphi$, as well as of the solution $\Psi$, and from equations~\eqref{eq:GL} and~\eqref{eq:GP} that the function $\varepsilon(x, t) := e^{- i \varphi(t)} \Psi(x + a(t), t) - V_1(x)$ in~\eqref{def:eps} satisfies
\begin{equation}
\label{eq:eps}
i \partial_t \varepsilon + \Delta \varepsilon + \big( 1 - |V_1|^2 \big) \varepsilon + \eta_\varepsilon \big( V_1 + \varepsilon) - \varphi'(t) \big( V_1 + \varepsilon \big) - i a'(t) \cdot \big( \nabla V_1 + \nabla \varepsilon \big) = 0,
\end{equation}
where $\eta_\varepsilon = 1 - |V_1 + \varepsilon|^2 - (1 - |V_1|^2)$, as before. With this equation at hand, we can differentiate the three orthogonality conditions in~\eqref{eq:orthocond} in order to obtain the system
\begin{equation}
\label{eq:sys-modul-deriv}
\boM_\varepsilon(t) \begin{pmatrix} a_1'(t) \\ a_2'(t) \\ \varphi'(t) \end{pmatrix} = \boF_\varepsilon(t).
\end{equation}
Since
$$
\int_{\R^2} \chi \langle \partial_x V_1, \partial_y V_1 \rangle_\C = \int_{\R^2} \chi \langle i V_1, \partial_x V_1 \rangle_\C = \int_{\R^2} \chi \langle i V_1, \partial_y V_1 \rangle_\C = 0,
$$
the matrix $\boM_\varepsilon$ in this formula is given by
\begin{equation}
\label{def:boM-eps}
\boM_\varepsilon = \boM_0 + \begin{pmatrix} \int_{\R^2} \chi \langle \partial_x V_1, \partial_x \varepsilon \rangle_\C & \int_{\R^2} \chi \langle \partial_x V_1, \partial_y \varepsilon \rangle_\C & \int_{\R^2} \chi \langle i \partial_x V_1, \varepsilon \rangle_\C \\ \int_{\R^2} \chi \langle \partial_y V_1, \partial_x \varepsilon \rangle_\C & \int_{\R^2} \chi \langle \partial_y V_1, \partial_y \varepsilon \rangle_\C & \int_{\R^2} \chi \langle i \partial_y V_1, \varepsilon \rangle_\C \\ \int_{\R^2} \chi \langle V_1, i \partial_x \varepsilon \rangle_\C & \int_{\R^2} \chi \langle V_1, i \partial_y \varepsilon \rangle_\C & \int_{\R^2} \chi \langle V_1, \varepsilon \rangle_\C \end{pmatrix},
\end{equation}
with
\begin{equation}
\label{def:boM-0}
\boM_0 = \begin{pmatrix} \int_{\R^2} \chi \big|\partial_x V_1|^2 & 0 & 0 \\ 0 & \int_{\R^2} \chi |\partial_y V_1|^2 & 0 \\ 0 & 0 & \int_{\R^2} \chi |V_1|^2 \end{pmatrix}.
\end{equation}
Similarly, the right-hand side $\boF_\varepsilon$ is equal to
\begin{equation}
\label{def:boF-eps}
\boF_\varepsilon = \begin{pmatrix} \int_{\R^2} \Big( \chi \big( (1 - |V_1|^2) 
\langle i \partial_x V_1, \varepsilon \rangle_\C + \eta_\varepsilon \langle i 
\partial_x V_1, V_1 + \varepsilon \rangle_\C \big) - \langle i \nabla \big( \chi 
\partial_x V_1 \big), \nabla \varepsilon \rangle_\C \Big) \\ \int_{\R^2} \Big( 
\chi \big( (1 - |V_1|^2) \langle i \partial_y V_1, \varepsilon \rangle_\C + 
\eta_\varepsilon \langle i \partial_y V_1, V_1 + \varepsilon \rangle_\C \big) - 
\langle i \nabla \big( \chi \partial_y V_1 \big), \nabla \varepsilon \rangle_\C 
\Big) \\ \int_{\R^2} \Big( \chi \big( (1 - |V_1|^2) \langle V_1, \varepsilon 
\rangle_\C + \eta_\varepsilon \langle V_1, V_1 + \varepsilon \rangle_\C \big) - 
\langle \nabla \big( \chi V_1 \big), \nabla \varepsilon \rangle_\C \Big) 
\end{pmatrix}.
\end{equation}

At this stage, recall the existence of a constant $A > 0$ such that
$$
\big\| \varepsilon(\cdot, t) \big\|_H \leq A \alpha_\tau,
$$
by~\eqref{eq:est-modul}. Combining the definition of the norm $\| \cdot \|_H$, Lemma~\ref{lem:comp-emb} and the fact that the function $\chi$ is smooth and compactly supported, we deduce from~\eqref{def:boM-eps} the existence of $C > 0$ such that
\begin{equation}
\label{boule}
\big\| \boM_\varepsilon(t) - \boM_0 \big\| \leq C \alpha_\tau.
\end{equation}
Hence we can choose the value of $\tau$ small enough, so that the matrix $\boM_\varepsilon$ is invertible and its operator norm is less than a number $C > 0$ (depending only on $\alpha_\tau$).

Since $\eta_\varepsilon = - 2 \langle \varepsilon, V_1 \rangle_\C - |\varepsilon|^2$, we can similarly invoke the definition of the norm $\| \cdot \|_H$, Lemma~\ref{lem:comp-emb}, the fact that the function $\chi$ is smooth and compactly supported, as well as the Sobolev embedding theorem, in order to control the right-hand side $\boF_\varepsilon$ by
\begin{equation}
\label{et-bill}
\big| \boF_\varepsilon(t) \big| \leq C \| \varepsilon(\cdot, t) \|_H,
\end{equation}
for a further choice of the constant $C$. In view of~\eqref{eq:sys-modul-deriv} and~\eqref{boule}, we obtain
$$
|a'(t)| + |\varphi'(t)| \leq C \| \varepsilon(\cdot, t) \|_H.
$$
Estimate~\eqref{eq:est-evol-modul} then follows from the proofs of Theorems~\ref{thm:estimcore} and~\ref{thm:stab}.

We are now in position to complete the proof of 
Proposition~\ref{prop:evol-param} by a density argument. Consider an arbitrary 
initial datum $\Psi^0 \in \boE(\R)$ and the corresponding decomposition $\Psi_0 
= U_0 + w_0$ provided by Lemma~\ref{lem:pourCauchy}. We can find a sequence of 
functions $w_0^n$ in $H^3(\R)$ that tend to $w_0$ in $H^1(\R^2)$ as $n \to + 
\infty$. Moreover, the corresponding maps $a_n$ and $\varphi_n$ are in 
$\boC^1(\R, \R^2)$, respectively $\boC^1(\R, \R)$, and their time derivatives 
are given by~\eqref{eq:sys-modul-deriv}.

Recall here that the Gross-Pitaevskii flow is globally continuous with respect 
to the initial datum in $U_0 + H^1(\R)$ by Proposition~\ref{prop:cauchy2}. Going 
back to the proof of Lemma~\ref{lem:Hilbert-H}, we observe that the flow map 
$\Psi_0 \mapsto \Psi(\cdot, t)$ remains continuous from $U_0 + H^1(\R)$ to 
$\boC^0([- T, T], H)$ for any $T > 0$. In view of 
Proposition~\ref{prop:modul-param}, we first deduce that the maps $a_n$ and 
$\varphi_n$ converge in $\boC^0([- T, T], \R^2)$, respectively $\boC^0([- T, T], 
\R)$, to the maps $a$ and $\varphi$ corresponding to the initial datum $\Psi_0$. 
We also deduce that the function $\varepsilon$ in~\eqref{def:eps} also depends 
continuously in $\boC^0([- T, T], H)$ on the initial datum in $U_0 + H^1(\R)$. 
Going back to the proofs of~\eqref{boule} and~\eqref{et-bill}, it follows that 
the time derivatives $a_n'$ and $\varphi_n'$ are convergent in $\boC^0([- T, T], 
\R^2)$, respectively $\boC^0([- T, T], \R)$. This is enough to guarantee the 
continuously differentiability of the maps $a$ and $\varphi$, and that their 
time derivatives $a'$ and $\varphi'$ satisfy~\eqref{eq:sys-modul-deriv} in the 
limit $n \to + \infty$. Reproducing the proofs of~\eqref{boule} 
and~\eqref{et-bill}, we obtain~\eqref{eq:est-evol-modul} as before. This 
completes the proof of Proposition~\ref{prop:evol-param}. \qed

\appendix
\numberwithin{cor}{section}
\numberwithin{lem}{section}
\numberwithin{prop}{section}
\numberwithin{thm}{section}
\section{Properties of the vortex solution}
\label{sec:V-1}

Recall that the vortex solution $V_1$ takes the special form $V_1(x) = \rho_1(r) e^{i \theta}$ for any point $x = (r \cos(\theta), r \sin(\theta)) \in \R^2$. Several properties of its profile $\rho_1$ are useful in the course of our proofs. For the sake of completeness, we have collected 
them in the next lemma, as well as their consequences on the algebraic decay rate of the lower order derivatives of $V_1$. 

\begin{lem}
\label{lem:propprofil}
$(i)$ There exists a unique solution $\rho_1 : [0, \infty) \to \R$ to the 
ordinary differential equation
\begin{equation}
\label{eq:rho1}
\rho_1''(r) + \frac{\rho_1'(r)}{r} - \frac{\rho_1(r)}{r^2} + \rho_1(r) \big( 1 - \rho_1(r)^2 \big) = 0,
\end{equation}
with $\rho_1(0) = 0$, and $\rho_1(r) \to 1$ as $r \to + \infty$. The function $\rho_1$ is smooth, increasing and it satisfies
$$
\rho_1(r) = A_1 \Big( r - \frac{r^3}{8} + \boO(r^5) \Big) \text{ as } r 
\to 0,
$$
with $A_1 = \rho_1'(0) > 0$.

\noindent $(ii)$ Moreover, the function $\rho_1$ satisfies the asymptotics
$$
\rho_1(r) = 1 - \frac{1}{2 r^2} - \frac{9}{8 r^4} + \boO \Big( \frac{1}{r^6} \Big) ,
$$
as well as
$$
\rho_1'(r) = \frac{1}{r^3} + \frac{9}{2 r^5} + \boO \Big( \frac{1}{r^7} 
\Big), \quad \rho_1''(r) = - \frac{3}{r^4} - \frac{45}{2 r^6} + \boO \Big( \frac{1}{r^8} \Big) \quad \text{and} \quad \rho_1'''(r) \sim \frac{12}{r^5},
$$
in the limit $r \to + \infty$. In particular, we have
$$
1 - \rho_1(r)^2 = \frac{1}{r^2} + \frac{2}{r^4} + \boO \Big( \frac{1}{r^6} \Big),
$$
as $r \to + \infty$.

\noindent $(iii)$ As a consequence, there exists a universal constant $C > 0$ such that
$$
\big| \nabla V_1(x) \big| \leq \frac{C}{1 + |x|}, \quad \big| d^2 V_1(x) \big| \leq \frac{C}{1 + |x|^2} \quad \text{and} \quad \big| d^3 V_1(x) \big| \leq \frac{C}{1 + |x|^3},
$$
for any $x \in \R^2$.
\end{lem}

\begin{proof}
Statement $(i)$ is proved in~\cite{CheElQi1, HervHer1}. Statement $(ii)$ is given in~\cite[Theorem 3.4]{CheElQi1}, except the expansion of the function $1 - \rho_1^2$ that is a direct consequence of the one for $\rho_1$, and the asymptotics for the third derivative $\rho_1'''$ that is obtained by differentiating~\eqref{eq:rho1}.

Since the vortex solution $V_1$ is smooth on $\R^2$, the proof of Statement $(i)$ reduces to quantify the algebraic decay of the lower order derivatives of $V_1$. Going back to~\eqref{eq:nabla-V1}, a direct computation provides the existence of a universal number such that
$$
\big| \nabla V_1(x) \big| \leq C \Big( \rho_1'(|x|) + \frac{\rho_1(|x|)}{|x|} \Big), \quad \big| d^2 V_1(x) \big| \leq C \Big( |\rho_1''(|x|)| + \frac{\rho_1'(|x|)}{|x|} + \frac{\rho_1(|x|)}{|x|^2} \Big),
$$
and
$$
\big| d^3 V_1(x) \big| \leq C \Big( |\rho_1'''(|x|)| + \frac{|\rho_1''(|x|)|}{|x|} + \frac{\rho_1'(|x|)}{|x|^2} + \frac{\rho_1(|x|)}{|x|^3} \Big),
$$
for any $x \in \R^2$. The bounds in Statement $(iii)$ then follow from Statement $(ii)$. This concludes the proof of Lemma~\ref{lem:propprofil}.
\end{proof}

\section{Properties of the function spaces}

In this appendix, we gather some properties related to the Hilbert space $H$ and the metric space $E$ that are useful in the course of our proofs.

\subsection{Properties of the Hilbert space \textit{H}}

Recall that the vector space $H$ is defined as
$$
H := \Big\{ \psi \in L_\text{loc}^2(\R^2)\text{ s.t. }\nabla (\psi \bar{V}_1) \in L^2(\R^2) \text{ and } (1 - |V_1|^2)^\frac{1}{2} \, \nabla \psi \in 
L^2(\R^2) \Big\},
$$
and that it is naturally endowed with the scalar product
$$
\big\langle \psi_1, \psi_2 \big\rangle_H := \int_{\R^2} \Big( \langle \nabla (\psi_1 \bar{V}_1), \nabla (\psi_2 \bar{V}_1) \rangle_\C + (1 - |V_1|^2) \langle \nabla \psi_1, \nabla \psi_2 \rangle_\C \Big).
$$
Concerning its topological properties, we show

\begin{lem}
\label{lem:Hilbert-H}
The vector space $H$ is a Hilbert space for the scalar product $\langle \cdot, \cdot \rangle_H$. Moreover,
$$
H^1(\R^2) \subset H.
$$
\end{lem}

\begin{proof}
Consider a Cauchy sequence $(\psi_n)_{n \in \N}$ of $H$. Fix $R > 0$ and denote by
$$
m_R(f) = \frac{1}{\pi R^2} \int_{B_R} f,
$$
the average mean of a locally integrable function $f$ on the ball $B_R$. As a consequence of the Poincar\'e-Wirtinger inequality, the sequences $(\psi_n \bar{V}_1 - m_R(\psi_n \bar{V}_1))_{n \in \N}$ and $(\psi_n - m_R(\psi_n))_{n \in \N}$ are Cauchy sequences, therefore convergent sequences, in $L^2(B_R)$. Hence, the sequence $(m_R(\psi_n) \bar{V}_1 - m_R(\psi_n \bar{V}_1))_{n \in \N}$ is also convergent in this space. Taking 
the scalar product in $L^2(B_R)$ with the function $\bar{V}_1$, we deduce that $(m_R(\psi_n))_{n \in \N}$ is a convergent sequence, so that $(\psi_n)_{n \in \N}$ is convergent in $L^2(B_R)$, and more generally in $L_\text{loc}^2(\R^2)$.

Let us denote by $\psi_\infty$ its limit. Going back to the definition of 
a Cauchy sequence in $H$, we know that the sequences $(\nabla (\psi_n \bar{V}_1))_{n \in \N}$ and $((1 - |V_1|^2)^{1/2} \nabla \psi_n)_{n \in \N}$ 
are convergent in $L^2(\R^2)$. Testing them in front of smooth compactly supported functions and taking the limit $n \to + \infty$ in the sense of distributions, we check that their limits in $L^2(\R^2)$ are equal to $\nabla (\psi_\infty \bar{V}_1)$, respectively $(1 - |V_1|^2)^{1/2} \nabla \psi_\infty$. Hence the function $\psi_\infty$ lies in $H$ and it is the limit of the Cauchy sequence $(\psi_n)_{n \in \N}$ in this space.

When $\psi \in H^1(\R^2)$, we also have
$$
\nabla \big( \psi \bar{V}_1 \big) = \bar{V}_1 \nabla \psi + \psi \nabla 
\bar{V}_1 \in L^2(\R^2),
$$
since $V_1$ and $\nabla V_1$ are bounded functions. For the same reason, the function $(1- |V_1|^2)^{1/2} \nabla \psi$ is in $L^2(\R^2)$. This proves that $\psi$ belongs to $H$ and concludes the proof of Lemma~\ref{lem:Hilbert-H}.
\end{proof}

A drawback of the previous definition for the Hilbert space $H$ lies in the property that its canonical norm $\| \psi \|_H$ does not provide any direct control on the function $\psi$ under consideration, but only on its 
gradient. In order to recover such a control, we next establish that the Hilbert space $H$ compactly embeds into suitable weighted Lebesgue spaces. As a consequence of this result, we especially gain a local control on the functions in $H$ that turns out to be very useful in our proofs.

\begin{lem}
\label{lem:comp-emb}
Let $s > 1$. The Hilbert space $H$ continuously embeds into the weighted Lebesgue space
$$
L_{- s}^2(\R^2) := \Big\{ \psi \in L_\text{loc}^2(\R^2)\text{ s.t. }\| \psi \|_{L_{- s}^2}^2 := \int_{\R^2} \frac{|\psi(x)|^2}{(1+|x|^2)^s} \, dx< \infty \Big\}.
$$
Moreover this embedding is compact.
\end{lem}

\begin{proof}
Let $\psi \in H$. In view of Lemma~\ref{lem:propprofil}, there exists a universal constant $K > 0$ such that
$$
\int_{B_2} |\psi |^2 \leq K \int_{B_2} | \nabla \bar{V}_1|^2 \, 
|\psi|^2.
$$
Since $\psi \nabla \bar{V}_1 = \nabla (\bar{V}_1 \psi) - \bar{V}_1 \nabla \psi$, we obtain
$$
\int_{B_2} |\psi |^2 \leq 2 K \int_{B_2} \Big( |\nabla (\bar{V}_1 \psi) |^2 + |V_1 |^2 |\nabla \psi|^2 \Big),
$$
Invoking once again Lemma~\ref{lem:propprofil}, we can find a further universal constant $K > 0$ such that
\begin{equation}
\label{Hardy0}
\int_{B_2} |\psi |^2 \leq K \int_{B_2} \Big( |\nabla (\bar{V}_1 
\psi)|^2 + (1 - | V_1 |^2) |\nabla \psi|^2 \Big) \leq K \| \psi \|_H^2.
\end{equation}
In particular, there exists a universal constant $K > 0$ such that
\begin{equation}
\label{Hardy3}
\| \nabla (\bar{V}_1 \psi) \|_{L^2(\R^2)} + \| \bar{V}_1 \psi \|_{L^2
 (B_2)} \leq K_0 \| \psi \|_H.
\end{equation}

Let us now show the existence of $K_s > 0$, depending only on $s$, such that
\begin{equation}
\label{Hardy2}
\int_{\R^2} \frac{|f(x)|^2}{(1 + | x |^2)^s} \, dx \leq K_s \big( \| \nabla f \|_{L^2(\R^2)}^2 + \| f \|_{L^2(B_2)}^2 \big)
\end{equation}
for any function $f \in \dot{H}^1 (\R^2)$. Consider, as previously, a smooth, decreasing cut-off function $\chi : \R_+ \to [0, 1]$ that satisfies $\chi \equiv 1$ on $[0, 1]$ and $\chi \equiv 0$ for $r \geq 2$. We first show the existence of $K_s > 0$ such that
\begin{equation}
\label{Hardy}
\int_{\R^2} \frac{| f (x) |^2}{| x |^{2 s}} (1 - \chi (| x |)) \, dx \leq 
K_s \big( \| \nabla f \|_{L^2 (\R^2)}^2 + \| f
 \|_{L^2 (B_2)}^2 \big).
\end{equation}
Indeed, we can combine the use of polar coordinates and an integration by 
parts in order to compute
$$
\int_{\R^2} \frac{| f (x) |^2}{| x |^{2 s}} (1 - \chi (| x |)) \, dx = - \frac{1}{2 (s - 1)} \int_0^{2 \pi} \int_0^{+ \infty} \frac{1}{r^{2 s - 2}} \partial_r \big( | f |^2 (1 - \chi (r)) \big) \, dr \, d\theta.
$$
We deduce that
\begin{align*}
\int_{\R^2} \frac{| f (x) |^2}{| x |^{2 s}} & (1 - \chi (| x |)) \, dx
\\
& \leq \frac{1}{2 (s - 1)} \int_{\R^2} \Big( \frac{|f(x)|^2}{| x |^{2 s - 1}} | 
\chi' (| x |) | + \frac{2 | f(x)| | \nabla f(x) |}{| x |^{2 s - 1}} (1 - \chi (| 
x |)) \Big) \, dx.
\end{align*}
Estimate~\eqref{Hardy} then follows from the Cauchy-Schwarz inequality that provides
$$
\int_{\R^2} \frac{| f(x) | | \nabla f(x) |}{| x |^{2 s - 1}} (1 - \chi(| x |)) \, dx \leq \bigg( \int_{\R^2} \frac{| f (x) |^2}{| x
 |^{2 s}} (1 - \chi (| x |)) \, dx \bigg)^\frac{1}{2} \| \nabla f \|_{L^2 
(\R^2)},
 $$
when $4 s - 2 > 2 s$. In turn, estimate~\eqref{Hardy2} follows from checking that
$$
\int_{\R^2} \frac{| f(x) |^2}{(1 + | x |^2)^s} \, dx \leq K \bigg( \| f \|_{L^2 (B_2)}^2 + \int_{\R^2} \frac{| f |^2}{| x
 |^{2 s}} (1 - \chi (| x |)) \, dx \bigg),
$$
for a further $K > 0$. In view of~\eqref{Hardy3}, we finally deduce from~\eqref{Hardy2} for $f = \bar{V}_1
\psi$ and from~\eqref{Hardy0} that the space $H$ continuously embeds into 
$L_{- s}^2
(\R^2)$. Let us now show that this embedding is compact.

Consider a bounded sequence $(\psi_n)_{n \in \N}$ of $H$. Given any $R > 0$, it follows from the previous embedding that this sequence 
is also bounded in $H^1(B_R)$. Invoking the Rellich-Kondrachov theorem and performing a diagonal argument, we can construct a function $\psi \in H_\text{loc}^1 (\R^2)$ such that, up to a subsequence,
\begin{equation}
\label{laurel}
\| \psi_n - \psi \|_{L^2 (B_R)} \to 0,
\end{equation}
when $n \to + \infty$, for any $R > 0$. By weak convergence, the function $\psi$ also belongs to $H$. In particular, since $(1 + s)/2 > 1$ for $s > 1$, we can invoke the previous continuous embedding in order to find $A > 0$ such that
$$
\int_{\R^2} \frac{|\psi_n(x)|^2}{(1 + | x |^2)^{\frac{1 + s}{2}}} \, dx + 
\int_{\R^2} \frac{|\psi(x)|^2}{(1 + | x |^2)^{\frac{1 + s}{2}}} \, dx \leq A,
$$
for any $n \in \N$. Now, take any $\delta > 0$, and let us conclude the proof of the compactness by showing the existence of a number $N_\delta \in \N$ such that
\begin{equation}
\label{HardyCompactness}
\int_{\R^2} \frac{|\psi_n(x) - \psi(x)|^2}{(1 + | x |^2)^s} \, dx \leq \delta,
\end{equation}
for any $n \geq N_\delta$. Indeed, we have
\begin{align*}
\int_{\R^2 \setminus B_R} \frac{|\psi_n(x) - \psi(x)|^2}{(1 + |x|^2)^s} \, dx \leq \frac{2}{R^\frac{s - 1}{2}} \int_{\R^2
 \setminus B_R} \frac{|\psi_n(x)|^2 + |\psi(x)|^2}{(1 + | x |^2)^{\frac{1 + s}{2}}} \, dx \leq \frac{2 A}{R^\frac{s - 1}{2}},
\end{align*}
for any $R > 0$. Since $(s - 1)/2 > 0$, we can choose $R_\delta$ such that $2
A/R_\delta^{(s - 1)/2} \leq \delta/2$. Furthermore, we compute
$$
\int_{B_{R_\delta}} \frac{|\psi_n(x) - \psi(x)|^2}{(1 + | x |^2)^s} \, 
dx \leq \int_{B_{R_\delta}} |\psi_n - \psi|^2.
$$
Hence, we deduce from~\eqref{laurel} the existence of an integer $N_\delta$ such that
$$
\int_{B_{R_\delta}} |\psi_n - \psi|^2 \leq \frac{\delta}{2},
$$
for any $n \geq N$. This concludes the proof of~\eqref{HardyCompactness} and of Lemma~\ref{lem:comp-emb}.
\end{proof}

Another crucial property of the Hilbert space $H$ lies in the fact that it is 
left invariant by translations and phase shifts. More precisely, we can show the 
following estimates for these two classes of operations.

\begin{lem}
\label{lem:bdd-translation}
Let $d \in \R^2$ and $\vartheta \in \R$. Given any function $\psi \in H$, the 
functions $e^{- i \vartheta} \psi$ and $\psi(\cdot + d)$ also belong to $H$. 
Moreover,
$$
\big\| e^{- i \vartheta} \psi \big\|_H = \big\| \psi \big\|_H,
$$
and there exists $C(|d|) > 0$, depending continuously on the norm $|d|$ in 
$\R_+$, such that
$$
\big\| \psi(\cdot + d) \big\|_H \leq C(|d|) \big\| \psi \big\|_H.
$$
\end{lem}

\begin{proof}
Let $\psi \in H$. By definition, we have
$$
\big\| e^{- i \vartheta} \psi \big\|_H = \big\| \psi \big\|_H,
$$
so that the function $e^{- i \vartheta} \psi$ is also in $H$. Similarly, we will prove that the function $\psi(\cdot + d)$ belongs to $H$ by bounding the quantity
$$
\big\| \psi(\cdot + d) \big\|_H^2 = \int_{\R^2} \Big( \big| \nabla( \bar{V}_1(\cdot - d) \psi) \big|^2 + \big( 1 - |V_1(\cdot - d)|^2 \big) \big| \nabla \psi \big|^2 \Big) =: I(\psi).
$$
In this direction, we will split the integral $I(\psi)$ as $I(\psi) = I_1(\psi) + I_2(\psi) + I_3(\psi)$, where
$$
I_1(\psi) := \int_{B_{|d| + 1}} \big| \nabla( \bar{V}_1(\cdot - d) \psi) \big|^2, \quad I_2(\psi) := \int_{B_{|d| + 1}^c} \big| \nabla( \bar{V}_1(\cdot - d) \psi) \big|^2,
$$
and
$$
I_3(\psi) = \int_{\R^2} \big( 1 - |V_1(\cdot - d)|^2 \big) \big| \nabla 
\psi \big|^2.
$$

Concerning the integral $I_1(\psi)$, we check that
\begin{align*}
I_1(\psi) \leq & 2 \int_{B_{|d| + 1}} \Big( |\nabla V_1(\cdot - d)|^2 |\psi|^2 + |V_1(\cdot - d)|^2 |\nabla \psi|^2 \Big)\\
\leq & 2 \| \nabla V_1 \|_{L^\infty}^2 \int_{B_{|d| + 1}} |\psi|^2 + \frac{1}{1 - \rho_1(|d| + 1)^2} \int_{B_{|d| + 1}} (1 - |V_1|^2) |\nabla \psi|^2.
\end{align*}
In view of the proof of Lemma~\ref{lem:comp-emb}, there exists a positive 
number $C(|d|)$, depending continuously on $|d|$ in $\R_+$, such that
$$
\int_{B_{|d| + 1}} |\psi|^2 \leq C(|d|) \| \psi \|_H^2,
$$
so that
\begin{equation}
\label{eq:est-I1}
I_1(\psi) \leq \Big( 2 C(|d|) \| \nabla V_1 \|_{L^\infty}^2 + \frac{1}{1 - \rho_1(|d| + 1)^2} \Big) \| \psi \|_H^2.
\end{equation}

In order to estimate the integral $I_2(\psi)$, we use the property that the function $V_1$ only vanishes at the origin. We obtain
\begin{align*}
I_2(\psi) & = \int_{B_{|d| + 1}^c} \Big| \nabla \Big( \frac{\bar{V}_1(\cdot - d)}{\bar{V}_1} \bar{V}_1 \psi \Big) \Big|^2\\
& \leq 2 \int_{B_{|d| + 1}^c} \bigg( \Big| \nabla \Big( \frac{V_1(\cdot 
- d)}{V_1} \Big) \Big|^2 |V_1|^2 |\psi|^2 + \frac{|V_1(\cdot - d)|^2}{|V_1|^2} \big| \nabla (\bar{V}_1 \psi) \big|^2 \bigg).
\end{align*}
The second term in the last integral of the previous inequality is bounded by
\begin{equation}
\label{eq:est-I2a}
\int_{B_{|d| + 1}^c} \frac{|V_1(\cdot - d)|^2}{|V_1|^2} \big| \nabla (\bar{V}_1 \psi) \big|^2 \leq \frac{1}{\rho_1(|d| + 1)^2} \| \psi \|_H^2.
\end{equation}
Concerning the first term, we compute
$$
|V_1|^2 \Big| \nabla \Big( \frac{V_1(\cdot - d)}{V_1} \Big) \Big|^2 = \frac{1}{|V_1|^2} \Big| V_1 \nabla V_1(\cdot - d) -V_1(\cdot - d) \nabla V_1 \Big|^2,
$$
so that by~\eqref{eq:nabla-V1}, we obtain
\begin{align*}
|V_1(x)|^2 \Big| \nabla \Big( \frac{V_1(x - d)}{V_1(x)} \Big) \Big|^2 = 
& \rho_1(x - d)^2 \, \frac{|d|^2}{|x|^2 |x - d|^2} + \rho_1'(x - d)^2 + \frac{\rho_1(x - d)^2}{\rho_1(x)^2} \rho_1'(x)^2\\
& - 2 \frac{\rho_1(x - d) \, x \cdot (x - d)}{\rho_1(x) |x| |x - d|} \rho_1'(x) \rho_1'(x - d).
\end{align*}
We check that
$$
\frac{1}{|x - d|} \leq \frac{|d| + 1}{|x|},
$$
when $|x| \geq |d| + 1$. In view of Lemma~\ref{lem:propprofil}, we deduce 
that
\begin{equation}
\label{cros}
|V_1(x)|^2 \Big| \nabla \Big( \frac{V_1(x - d)}{V_1(x)} \Big) \Big|^2 \leq C \Big( \frac{|d|^2 (|d| + 1)^2}{|x|^4} + \frac{(|d| + 1)^6}{|x|^6} + \frac{1}{\rho_1(|d| + 1)^2 |x|^6} \Big),
\end{equation}
for some $C > 0$, not depending on $d$. Going back to the proof of 
Lemma~\ref{lem:comp-emb}, we infer the existence of a further $C(|d|) > 0$, 
depending continuously on $|d|$ in $\R_+$, such that
\begin{equation}
\label{eq:est-I2b}
\int_{B_{|d| + 1}^c} \Big| \nabla \Big( \frac{V_1(\cdot - d)}{V_1} \Big) \Big|^2 |V_1|^2 |\psi|^2 \leq C(|d|) \big\| \psi \big\|_H^2.
\end{equation}

Similarly, we bound the integral $I_3(\psi)$ by
$$
I_3(\psi) \leq \Big\| \frac{1 - |V_1(\cdot - d)|^2}{1 - |V_1|^2} \Big\|_{L^\infty(\R^2)} \int_{\R^2} (1 - |V_1|^2) \big| \nabla \psi \big|^2.
$$
When $|x| \leq |d| + 1$, we observe that
$$
\Big| \frac{1 - |V_1(x - d)|^2}{1 - |V_1(x)|^2} \Big| \leq \frac{1}{1 - \rho_1(|d| + 1)^2}.
$$
On the other hand, we deduce from Lemma~\ref{lem:propprofil} the existence of $C > 0$, 
not depending on $d$, such that
$$
\Big| \frac{1 - |V_1(x - d)|^2}{1 - |V_1(x)|^2} \Big| \leq C \frac{|x|^2}{|x - d|^2} \leq 2 C \Big( 1 + \frac{|d|^2}{|x - d|^2}\Big) \leq 2 C (1 + 
|d|^2),
$$
for $|x| \geq 1 + |d|$. As a consequence, we obtain
$$
I_3(\psi) \leq \max \Big\{ \frac{1}{1 - \rho_1(|d| + 1)^2}, 2 C (1 + |d|^2) \Big\} \big\| \psi \big\|_H^2.
$$
The conclusion then follows from~\eqref{eq:est-I1},~\eqref{eq:est-I2a} and~\eqref{eq:est-I2b}.
\end{proof}

Concerning translation, we can refine the estimate in Lemma~\ref{lem:bdd-translation} in the 
special case of the vortex solution $V_1$. In the next lemma, we establish some local Lipschitz 
continuity of the function $V_1(\cdot + d)$ with respect to the translation parameter $d$. 
This property is useful in the previous construction of the modulation parameters.

\begin{lem}
\label{lem:Lipschitz-translation}
Let $d \in \R^2$. There exists $C(|d|) > 0$, depending continuously on the norm 
$|d|$ in $\R_+$, such that
$$
\big\| V_1(\cdot + d) - V_1 \big\|_H \leq C(|d|) \, |d|.
$$
\end{lem}

\begin{proof}
By definition, we have
$$
\big\| V_1(\cdot + d) - V_1 \big\|_H^2 = \int_{\R^2} \Phi_{V_1}(x, d) \, dx =: 
I_{V_1}(d),
$$
where
$$
\Phi_{V_1}(x, d) := \big| \nabla \big( \bar{V}_1(x) (V_1(x +d) - V_1(x)) \big) 
\big|^2 + (1 - |V_1(x)|^2) \big| \nabla (V_1(x +d) - V_1(x)) \big|^2,
$$
for any pair $(x, d) \in \R^2 \times \R^2$. The proof is then based on the 
property that the integral $I_{V_1}$ is of class $\boC^2$ on $\R^2$, with 
$I_{V_1}(0) = 0$ and $\nabla I_{V_1}(0) = 0$. In this case, we can apply the 
Taylor formula in order to obtain
$$
\big\| V_1(\cdot + d) - V_1 \big\|_H^2 \leq M(|d|) |d|^2,
$$
with $M(|d|) := \max_{|x| \geq |d|} \| d^2 I_{V_1}(x) \|$ being continuous with 
respect to $|d|$. It is then enough to take the square root of this inequality 
in order to complete the proof of Lemma~\ref{lem:Lipschitz-translation}.

Hence, we are reduced to check the second order continuous differentiability of 
the integral $I_{V_1}$ by applying the dominated convergence theorem. In view of 
Lemma~\ref{lem:propprofil}, the function $\Phi_{V_1}$ is smooth on $\R^2 \times 
\R^2$, with
\begin{equation}
\label{tolofua}
\begin{split}
\partial_{d_i} \Phi_{V_1}(x, d) := & 2 \big\langle \nabla \big( \bar{V}_1 (V_1(x 
+d) - V_1(x)) \big), \nabla \big( \bar{V}_1(x) \partial_{x_i} V_1(x +d) \big) 
\big\rangle_\C\\
& + 2 (1 - |V_1(x)|^2) \big\langle \nabla (V_1(x +d) - V_1(x)), \partial_{x_i} 
\nabla V_1(x + d) \big\rangle_\C,
\end{split}
\end{equation}
and
\begin{equation}
\label{woki}
\begin{split}
\partial_{d_i} \partial_{d_j} \Phi_{V_1}(x, d) := & 2 \big\langle \nabla \big( 
\bar{V}_1 \partial_{x_j} V_1(x +d) \big), \nabla \big( \bar{V}_1(x) 
\partial_{x_i} V_1(x +d) \big) \big\rangle_\C\\
& + 2 \big\langle \nabla \big( \bar{V}_1 (V_1(x +d) - V_1(x)) \big), \nabla 
\big( \bar{V}_1(x) \partial_{x_i} \partial_{x_j} V_1(x +d) \big) 
\big\rangle_\C\\
& + 2 (1 - |V_1(x)|^2) \big\langle \partial_{x_j} \nabla V_1(x +d), 
\partial_{x_i} \nabla V_1(x + d) \big\rangle_\C\\
& + 2 (1 - |V_1(x)|^2) \big\langle \nabla (V_1(x +d) - V_1(x)), \partial_{x_i} 
\partial_{x_j} \nabla V_1(x + d) \big\rangle_\C,
\end{split}
\end{equation}
for $1 \leq i, j \leq 2$ and $(x, d) \in \R^2 \times \R^2$. Recall that
$$
\frac{1}{|x + d|} \leq \frac{|d| + 1}{|x|},
$$
when $|x| \geq |d| + 1$. In view of Lemma~\ref{lem:propprofil}, this inequality 
is enough to find $C(|d|) > 0$, depending continuously on $|d|$ in $\R_+$, such 
that
\begin{align*}
\big| V_1(x + d) \big| & + \big( 1 + |x|) \big| \nabla V_1(x + d) \big| \\
& + \big( 1 + |x|^2) \big| d^2 V_1(x + d) \big| + \big( 1 + |x|^3) \big| d^3 
V_1(x + d) \big| \leq C(|d|).
\end{align*}
We then infer from~\eqref{tolofua} and~\eqref{woki} the bounds
$$
\big| \partial_{d_i} \Phi_{V_1}(x, d) \big| + (1 + |x|) \big| \partial_{d_i} 
\partial_{d_j} \Phi_{V_1}(x, d) \big| \leq \frac{C(|d|)}{1 + |x|^3}.
$$
Arguing as for the proof of~\eqref{cros}, we also obtain the refined bound
$$
\big| \Phi_{V_1}(x, d) \big| \leq \frac{C(|d|)}{1 + |x|^4}.
$$
Applying the dominated convergence theorem with $d$ lying in bounded subsets of 
$\R^2$, we conclude that the integral $I_{V_1}$ is of class $\boC^2$ on $\R^2$. 
The facts that $I_{V_1}(0) = 0$ and $\nabla I_{V_1}(0) = 0$ follows from the 
identities $\Phi_{V_1}(x, 0) = \partial_{d_1} \Phi_{V_1}(x, 0) = \partial_{d_2} 
\Phi_{V_1}(x, 0) = 0$ for any $x \in \R^2$. This completes the proof of 
Lemma~\ref{lem:Lipschitz-translation}.
\end{proof}

\subsection{Properties of the metric space \textit{E}}

We now turn to the energy set
$$
E = \big\{ \psi \in H \text{ s.t. } 1 - |\psi|^2 \in L^2(\R^2) \big\},
$$
that we have endowed with the distance
$$
d_E \big( \psi_1, \psi_2 \big) = \big\| \psi_1 - \psi_2 \big\|_H + \big\| 
|\psi_2|^2 - |\psi_1|^2 \big\|_{L^2}.
$$
Recall that the classical Hamiltonian framework for the Gross-Pitaevskii 
equation is given by the set of functions with finite Ginzburg-Landau energy 
(see e.g.~\cite{BetGrSa2, Gerard2} and the references therein). The introduction 
of the energy set $E$ is reminiscent from this framework. Roughly speaking, this 
set is composed of functions with infinite Ginzburg-Landau energy due to a 
topological degree equal to $1$ at infinity, but of finite Ginzburg-Landau 
energy when this degree is suitably brought back to $0$. This interpretation can 
be made more effective through the following observation.

\begin{lem}
\label{lem:finiGL}
Let $\psi \in E$. The function $\psi \bar{V}_1$ has finite Ginzburg-Landau energy.
\end{lem}

\begin{proof}
By definition of the energy set $E$, the function $\nabla (\psi \bar{V}_1)$ is 
square integrable. For the potential term, we write
$$
1 - |\psi \bar{V}_1|^2 = |V_1|^2 (1 - |\psi|^2) + 1 - |V_1|^2.
$$
The right-hand side of this formula is also square integrable since the function 
$|V_1| = \rho_1$ is bounded by $1$ and the function $1 - |V_1|^2$ is square 
integrable. Hence, the Ginzburg-Landau energy of $\psi \bar{V}_1$ is finite.
\end{proof}

As a consequence of Lemma~\ref{lem:finiGL}, it is natural to rely on earlier 
results about the functions with finite Ginzburg-Landau energy (in particular 
in~\cite{Gerard2}) in order to describe the main properties of the energy set 
$E$. Our first result in this direction is

\begin{lem}
\label{lem:comp-E}
The energy set $E$ is a complete metric space for the distance $d_E$. Moreover, 
it satisfies
\begin{equation}
\label{eq:inclu-E}
E \subset L^4(\R^2) + L^\infty(\R^2) \subset L^2(\R^2) + L^\infty(\R^2),
\end{equation}
as well as
\begin{equation}
\label{eq:E+H1}
E + H^1(\R^2) = E.
\end{equation}
\end{lem}

\begin{proof}
Observe first that a Cauchy sequence $(\psi_n)_{n \in \N}$ of $E$ is a Cauchy 
sequence of the Hilbert space $H$. As a consequence of 
Lemmas~\ref{lem:Hilbert-H} and~\ref{lem:comp-emb}, it is convergent in $H$ and 
in $L_\text{loc}^2(\R^2)$ towards a limit function $\psi_\infty \in H$. Since 
$(1 - |\psi_n|^2)_{n \in \N}$ is a Cauchy sequence of $L^2(\R^2)$, it is also 
convergent in $L^2(\R^2)$. Moreover, its limit function is necessarily equal to 
$1 - |\psi_\infty|^2$ by almost everywhere convergence. Hence, the function 
$\psi_\infty$ is in $E$ and it is the limit of the sequence $(\psi_n)_{n \in 
\N}$ in this metric space that is therefore complete.

In order to prove~\eqref{eq:inclu-E}, we write a given function $\Psi \in E$ as 
$\Psi = \Psi \, \mathbbm{1}_{|\Psi| \geq 2} + \Psi \, \mathbbm{1}_{|\Psi| < 2}$. 
The second function in this decomposition is bounded. The first one is in 
$L^4(\R^2)$. Indeed, we know that $|\Psi| \leq 2 (|\Psi|^2 - 1)^{1/2}$ whenever 
$|\Psi| \geq 2$, and that $(|\Psi|^2 - 1)^{1/2} \in L^4(\R^2)$ by definition of 
$E$. This proves that $E \subset L^4(\R^2) + L^\infty(\R^2)$. The other 
inclusion in~\eqref{eq:inclu-E} then follows from the general property that 
$L^4(\R^2) + L^\infty(\R^2) \subset L^2(\R^2) + L^\infty(\R^2)$.

Concerning~\eqref{eq:E+H1}, Lemma~\ref{lem:Hilbert-H} guarantees that $H^1(\R^2) 
\subset H$. Therefore, we are reduced to establish that $1 - |\psi + u|^2$ is in 
$L^2(\R^2)$ when $\psi \in E$ and $u \in H^1(\R^2)$. For that purpose, we write
$$
|\psi + u|^2 - 1 = |\psi|^2 - 1 + |u|^2 + 2 \langle u, \psi \rangle_\C.
$$
In this formula, the function $1 - |\psi|^2$ is in $L^2(\R^2)$ by definition of 
the metric space $E$, so as the function $|u|^2$ due to the Sobolev embedding 
theorem. Since $u \in L^2(\R^2) \cap L^4(\R^2)$ again by the Sobolev embedding 
theorem, and $\psi \in L^\infty(\R^2) + L^4(\R^2)$ by~\eqref{eq:inclu-E}, the 
function $\langle u, \psi \rangle_\C$ also belongs to $L^2(\R^2)$. Hence, the 
function $|\psi + u|^2 - 1$ is in $L^2(\R^2)$, which completes the proof of 
Lemma~\ref{lem:comp-E}.
\end{proof}

In the spirit of Lemma~\ref{lem:Lipschitz-translation}, we also need at some
point the following Lipschitz estimate for translations on $V_1$ in $E.$

\begin{lem}
\label{lem:lipetaV1}
Let $d \in \R^2$. There exists $C(|d|) > 0$, depending continuously on the norm 
$|d|$ in $\R_+$, such that
$$
\big\| |V_1(\cdot + d)|^2 - |V_1|^2 \big\|_{L^2} \leq C(|d|) \, |d|.
$$
\end{lem}

\begin{proof}
The proof is very similar to the one of Lemma~\ref{lem:Lipschitz-translation}, 
and therefore we omit it.
\end{proof}

In order to tackle the Cauchy problem in the energy set $E$, we now relate it 
with the functional framework introduced in~\cite{BethSme1} to solve this 
problem for functions with non-zero degree at infinity. Recall that this 
framework was based on the set
\begin{equation}
\label{def:U}
\begin{split}
\boU(\R^2) := \big\{ U \in E \text{ s.t. } U \in L^\infty(\R^2), & \nabla |U| 
\in L^2(\R^2),\\ & \text{ and } \nabla^k U \in L^2(\R^2) \text{ for all } k \geq 
2 \big\}.
\end{split}
\end{equation}
Given a fixed function $U \in \boU(\R^2)$, it was proved in~\cite{BethSme1} that 
the Cauchy problem for~\eqref{eq:GP} is globally well-posed in $U + H^1(\R^2)$. 
This result can be applied in the context of the energy set $E$ since any 
function in this set can be decomposed as a function in $\boU(\R^2)$ plus an 
$H^1$-function.

\begin{lem}
\label{lem:pourCauchy}
Let $\psi \in E$. There exist two functions $U \in \boU(\R^2)$ and $w \in 
H^1(\R^2)$ such that
$$
\psi = U + w.
$$
\end{lem}

\begin{proof}
We fix a smooth, non-negative, compactly supported mollifier $\rho$. Given an 
arbitrary function $\psi \in E$, we decompose it as $\psi = U + w$, where
$$
U := V_1 \Big( 1 + \rho \ast \big( (\psi - V_1)\bar{V}_1 \big) \Big). 
$$
We first show that $w$ is in $H^1(\R^2)$. Setting $\varepsilon := \psi - V_1$, 
we compute
$$
w = \varepsilon - \big( \rho \ast (\varepsilon \bar{V}_1) \big) \, V_1 = \big( 
\varepsilon \bar{V}_1 - \rho \ast (\varepsilon \bar{V}_1) \big) V_1 + 
\varepsilon (1- |V_1|^2).
$$
Since $\int_{\R^2} \rho = 1$, we can find $C > 0$, depending only on $\rho$, 
such that
$$
\| f - \rho \ast f \|_{L^2} \leq C \| \nabla f \|_{L^2},
$$
for any function $f \in \dot{H}^1(\R^2)$. We use this inequality for $f = 
\varepsilon \bar{V}_1$, which belongs to $\dot{H}^1(\R^2)$ by definition of the 
vector space $H$. Combined with the facts that $V_1 \in L^\infty(\R^2)$, $1 - 
|V_1|^2 \in L^2(\R^2) \cap L^\infty(\R^2)$ and $\varepsilon \in L^2(\R^2) + 
L^\infty(\R^2)$ by Lemma~\ref{lem:comp-E}, we infer that $w$ is in $L^2(\R^2)$. 
In order to prove that $\nabla w \in L^2(\R^2)$, we use similar arguments, and 
the fact that $|\nabla \varepsilon| (1 - |V_1|^2) \leq |\nabla \varepsilon| (1 - 
|V_1|^2)^{1/2}$, where the latter term is in $L^2(\R^2)$ by definition of $H$.

We next show that $U \in \boU(\R^2)$. Since $|V_1| \leq 1$, we first have
$$
|U| \leq 1 + \rho \ast |\varepsilon|.
$$
Since $\varepsilon \in L^2(\R^2) + L^\infty(\R^2)$ by Lemma~\ref{lem:comp-E}, 
and $\rho \in L^\infty(\R^2) \cap L^2(\R^2)$, we deduce that $U \in 
L^\infty(\R^2)$. Similarly, we compute
$$
\big| \nabla |U| \big| \leq \big| \nabla |V_1| \big| \, \big( 1+\rho \ast 
|\varepsilon| \big) + \big| \rho \ast \nabla (\varepsilon \bar{V}_1) \big|. 
$$
By Lemma~\ref{lem:propprofil}, $\nabla |V_1|$ belongs to $L^2(\R^2)$, while 
$\nabla (\varepsilon \bar{V}_1)$ is in $L^2(\R^2)$ by definition of $H$. Hence, 
$\nabla |U|$ is also in $L^2(\R^2)$. The fact that $1-|U|^2 \in L^2(\R^2)$ is 
then a consequence of the fact that $\psi \in E$, $w \in H^1(\R^2)$ and 
Lemma~\ref{lem:comp-E}. It remains to show that $\nabla^k U \in L^2(\R^2)$ for 
any $k \geq 2$. In this direction, standard tame estimates yield
$$
\| \nabla^k U \|_{L^2} \leq C \Big( \| \nabla^k V_1 \|_{L^2} \big\| 1 + \rho 
\ast (\varepsilon \bar{V}_1) \big\|_{L^\infty} + \| V_1 \|_{L^\infty} \big\| 
\nabla ^{k - 1} \rho \ast \nabla(\varepsilon \bar{V}_1) \big\|_{L^2} \Big),
$$
and this quantity is finite by Lemma~\ref{lem:propprofil} and by definition of 
$H$. Finally, the function $U$ is also in $E$ by~\eqref{eq:E+H1}.
\end{proof}

As a direct consequence of Lemma~\ref{lem:pourCauchy}, we deduce that smooth 
functions are dense in the energy set $E$.

\begin{cor}
\label{cor:smooth-dense-E}
Let $\psi \in E$. There exist smooth functions $\psi_n \in E$ such that
$$
d_E(\psi_n, \psi) \to 0,
$$
as $n \to + \infty$.
\end{cor}

\begin{proof}
In view of Lemma~\ref{lem:pourCauchy}, we can find two functions $U \in 
\boU(\R^2) \cap E$ and $w \in H^1(\R^2)$ such that $\psi = U + w$. Since smooth, 
compactly supported functions are dense in $H^1(\R^2)$, there exist functions 
$w_n \in \boC_c^\infty(\R^2)$ such that $w_n \to w$ in $H^1(\R^2)$ as $n \to + 
\infty$. Set
$$
\psi_n = U + w_n,
$$
for any $n \in \N$. Since $U$ is smooth by definition of $\boU(\R^2)$, the 
functions $\psi_n$ are smooth. They also belong to $E$ by~\eqref{eq:E+H1}. 
Moreover, we have $\psi_n - \psi = w_n - w$. Arguing as in the proof of 
Lemma~\ref{lem:Hilbert-H}, we infer that
$$
\| \psi_n - \psi \|_H \leq C \| w_n - w \|_{H^1} \to 0,
$$
as $n \to + \infty$. Similarly, we have
$$
1 - |\psi_n|^2 - \big( 1 - |\psi|^2 \big) = 2 \langle U, w - w_n \rangle_\C + 
|w|^2 - |w_n|^2.
$$
Invoking the Sobolev embedding theorem is enough to prove the convergence of 
this quantity towards $0$ in $L^2(\R^2)$. This concludes the proof of 
Corollary~\ref{cor:smooth-dense-E}.
\end{proof}

\section{On the Cauchy problem in the space \textit{E}}
\label{sect:cauchy}

The goal of this section is to make a link with the analysis of the Cauchy 
problem for the Gross-Pitaevskii equation in~\cite{BethSme1}, and to prove 
Proposition~\ref{prop:cauchy}. We recall that in~\cite{BethSme1}, the Cauchy 
problem for~\eqref{eq:GP} is proved to be globally well-posed in $U + H^1(\R^2, 
\C)$ for any function $U \in \boU$, where the space $\boU$ is given 
in~\eqref{def:U}. Moreover, (some slightly different version of) the 
renormalized energy $\boE$ is shown to be preserved along the flow. In order to 
prove Proposition~\ref{prop:evol-param}, we need the following close extension 
of this result.

\begin{prop}
\label{prop:cauchy2}
Let $k \geq 1$ and $U \in \boU$. The Cauchy problem for~\eqref{eq:GP} is 
globally well-posed in $U + H^k(\R^2, \C)$, and the renormalized energy $\boE$ 
is conserved by the flow.
\end{prop}

\begin{proof}
Global well-posedness for $k = 1$ and $k = 2$ was proved in~\cite{BethSme1}. 
Concerning the renormalized energy, it is defined in~\cite{BethSme1} as
$$
\boE_U(w) := \frac{1}{2} \int_{\R^2} |\nabla w|^2 - \int_{\R^2} \langle \Delta U 
, w \rangle_\C + \frac{1}{4} \int_{\R^2} (1 - |U + w|^2)^2,
$$
for $\Psi = U + w \in U + H^1(\R^2, \C)$. Moreover, it was shown that
$$
\boE_U(w) = \lim_{r \to + \infty} \int_{B_r} \Big( e_\text{GL}(\Psi) - 
\frac{|\nabla U|^2}{2} \Big).
$$
In particular, we have
$$
\boE_U(w_t) = \boE(\Psi_t) - \boE(U) + \frac{1}{4} \int_{\R^2} (1 - |U|^2)^2,
$$
and therefore, $\boE(\Psi_t)$ is constant, since $U$ is fixed and $\boE_U$ is 
preserved by the flow.

Concerning the cases $k \geq 3$,~\footnote{We have made use of the case $k = 3$ 
in the course of the proof of Proposition~\ref{prop:evol-param}} local 
well-posedness follows as for $k = 2$, since the nonlinearity is Lipschitz due 
to the Sobolev embedding theorem of $H^k(\R^2)$ into $L^\infty(\R^2)$ when $k 
\geq 2$. Global existence for $k \geq 3$ is then a consequence of standard 
energy estimates and global existence for $k = 2$, again using the control of 
the uniform norm provided by the $H^2$-norm.
\end{proof}

Proposition~\ref{prop:cauchy} is finally a direct consequence of 
Lemma~\ref{lem:pourCauchy} and Proposition~\ref{prop:cauchy2}.

\begin{merci}
E.P. is supported by Tamkeen under the NYU Abu Dhabi Research Institute grant 
CG002. P.G. and D.S. acknowledge support from the project ``Dispersive and 
random waves'' (ANR-18-CE40-0020-01) of the Agence Nationale de la Recherche.
\end{merci}

\bibliographystyle{plain}
\bibliography{Bibliogr}

\end{document}